\documentstyle[12pt]{article}
\def\Br{{\rm Br}}
\def\Pic{{\rm Pic}}
\def\Pr{{\it Proof.~}}
\def\III{{\rm III}}

\def\halm{\vrule height1ex width.9ex depth.1ex}% halmos sign
\def\db #1 {{\bf#1}}

\def\Ker{{\rm Ker\;}}
\def\Im{{\rm Im}}
\def\Hom{{\rm Hom\,}}
\def\SL{{\rm SL}}

\def\A{{\rm A}}

\def\D{{\rm D}}
\def\E{{\rm E}}

\def\H{{\rm H}}

\def\K{{\rm K}}

\def\R{{\rm R}}

\font\germ=eufm10
\def\g{{\germ g}}
\def\h{{\germ h}}
\def\p{{\germ p}}
\def\gm{$\germ g\rm$}

\def\Coker{{\rm Coker}}
\def\Im{{\rm Im\;}}

\def\UK{{\rm UK}}

\def\H{{\rm H}}
\def\SL{{\rm SL}}
\def\SK{{\rm SK}}

\def\USK{{\rm USK}}
\def\K{{\rm K}}
\def\SU{{\rm SU}}
\begin{document}
\title{Weak approximation, Brauer and R-equivalence
in algebraic groups over arithmetical fields}
\author{Nguy\^e\~n Qu\^o\'{c} Th\v{a}\'ng\thanks{
Supported by a Lady Davis Fellowship and ICTP}}
\date{}
\maketitle
\begin{center}
Department of Mathematics, Israel Institute of Technology - Technion
\\
Haifa - 32000, Israel
\end{center}
\begin{abstract} We prove  
some  new relations between weak approximation
and some rational equivalence relations
(Brauer and R-equivalence) in algebraic groups
over arithmetical fields. By using weak approximation
and local - global approach, we compute completely
the group of Brauer equivalence classes of connected linear algebraic
groups over number fields, and also completely compute
the group of R-equivalence
classes of connected linear algebraic groups $G$, which either are defined
over a totally imaginary number field, or contains no anisotropic almost simple
factors of exceptional type $^{3,6}\D_4$, nor $\E_6$.
We discuss some consequences derived  from these, e.g., by giving
some new criteria for weak
approximation in algebraic groups over number fields,
by indicating a new way to give examples of non stably rational
algebraic groups over local fields and application to norm principle.
\end{abstract}
{\bf Introduction.} Let $G$ be a linear algebraic group defined
over a field $k$. There are two closely related questions in the arithmetic
theory of algebraic groups over fields :
the question of  weak approximation and that of rationality of a given $G$.
It is very difficult to study such questions for arbitrary groups
over arbitrary fields. One should restrict to some class of groups and
fields which are convenient in application.
\\
Let $X$ be a smooth algebraic variety defined over a field $k$,
$\bar X = X \times \bar k$, where $\bar k$ is an algebraic
closure of $k$. Denote by
$\Br(X)$ the cohomological Brauer group $\H^2_{et}(X,{\bf G}_m)$ of
$X$, $\Br_1(X) = \Ker (\Br(X) \to \Br(\bar X))$, $\Br_0 (X) = \Im (\Br (k)
\to \Br( X))$.
Following Manin, Colliot-Th\'el\`ene and Sansuc
(see [M1,2], [CTS1,2]), one defines the Brauer equivalence and
R-equivalence as follows. First one defines a pairing
$$X(k) \times \Br(X) \to \Br(k),
(x,b) \mapsto b(x),$$
where $b \in \H^2_{et}(X,{\bf G}_m)$ is considered (see [S], Section 6)
as an equivalence class of Azumaya algebras over $X$ and $b(x)$
is the equivalence class of central simple algebras over $k$, which is
considered as an element of $\Br(k)$.
\\
Two points $x,y \in X(k)$ are said to be {\it Brauer equivalent}
({\it Br-equivalent})
(resp. {\it R-equivalent}) if for any $b \in \Br(X)$, we have $b(x) = b(y)$
(resp. if there is a sequence 
of points $z_i \in X(k)$, $x=z_1, y=z_n$, such that 
for each pair $z_i,z_{i+1}$ there is
a $k$-rational map $f : {\bf P}^1 \to X$, regular at 0 and 1,
with $f(0)=z_i, f(1) = z_{i+1}, 1 \le i \le n-1$). In the definition
of Brauer equivalence above, one may also restrict to the subgroup $\Br_1 X$
of $\Br X$ to get a weaker equivalence relation. However, if $\bar X$ is 
rational over $k$ (which is the main case we are interested in),
it is known    
(cf. e.g. [CTS1], Lemme 16) 
that these two notions coincide. Moreover in [loc.cit],
Prop. 16, it was shown that the Brauer equivalence 
is weaker than R-equivalence, i.e.,
two points of $X(k)$, being R-equivalent, are necessarily Br-equivalent. 
In [loc.cit], basic theory of Brauer equivalence on 
tori defined over a field $k$ of
 characteristic 0 has been developed. In particular, in the arithmetic 
case, i.e.,
when $k$ is a local or global field, 
formulae for computations of the {\it group}
$T(k)/Br$ are given and
 it turns out to be in general a birational
invariant of $T$.
 Though $T(k)/Br$ is "computable",
 the group itself and its computation is in general non-trivial.
\\
In a subsequent paper [S],
 Sansuc developed  Brauer theory of linear algebaic groups $G$
over number fields, and applied it to obtain certain fundamental sequences
connecting various arithmetic (obstruction to weak approximation), 
cohomological (Tate - Shafarevich group) and geometric
 invariants (the first Galois
cohomology of the Picard group of a
smooth  compactification of $G$ over an algebraic closure
of $k$) for connected
linear algebraic groups $G$ over number fields $k$.
\\
 In this paper we continue the approach
taken by Colliot-Th\'el\`ene and Sansuc,
 to obtain certain connections between the above arithmetic, cohomological
and birational (Brauer  and R-) 
invariants of connected linear algebraic groups $G$ over local and
global fields of characteristic 0. As it was pointed out above,
in general, the group of
Brauer equivalence classes of $G$ is non-trivial, even in the case of tori. 
Therefore it is natural to ask what kind of analogs in the case of arbitrary
connected linear algebraic groups one can have. 
\\
Untill recently it was not known whether the group $G(k)/R$ is always
finite for any connected linear algebraic group $G$ defined over a number field
$k$. In his paper in 1993 ([G2]), Gille gave a proof of this 
finiteness property by using norm
and Hasse principles. However it was not known how one
can compute the actual group $G(k)/R$, nor even suggested how it might 
look like. 
The main result of this paper,
among others, is Theorem 4.12, which allows us not only to have a new
proof of this finiteness result for all connected linear algebaric
groups $G$, but also, by using the corresponding result for tori
done in [CTS1], to compute
completely the group $G(k)/R$
for those connected linear algebraic groups $G$,
which either  are defined over a totally
imaginary number field $k$, or contain no anisotropic almost simple factors
of exceptional type $^{3,6}\D_4$, nor $\E_6$. In these last 
critical cases, it reduces the computation
of $G(k)/R$ to a particular case of a well-known Platonov - Margulis
conjecture about the normal structure of almost simple simply connected
groups over number fields. 
\\
 Regarding the group $G(k)/Br$, we 
computed it completely, which in fact gives apriori (or preliminary)
information on the group $G(k)/R$.
The content of the paper after this introduction is as follows.
\\
\\
1. Recall of some basic facts from Brauer theory of algebraic groups.
\\
2. A Brauer relative of exact sequence (R) for algebraic tori.
\\
3. Some reductive analogs.
\\
4. Some variations and applications.
\\
5. An application to norms principle.
\\
6. Remarks, problems and conjectures.
\\
Bibliography.
\\
\\
{\it Notation.} 
Let $S$ be a finite set of valuations of a global field $k$, and $G$ a connected
linear algebraic group defined over $k$. Denote
by $Cl_S(G(k))$ the closure of $G(k)$ 
in the product topology, where $G(k)$ is embedded
diagonally into the direct product
$\prod_{v \in S}G(k_v)$
and $G(k_v)$ has the $v$-adic topology, induced from that of $k_v$.
We say that
$G$ has {\it weak approximation with respect to S} (or
{\it in S}) if $Cl_S(G(k))= \prod_{v \in S}G(k_v)$,
and {\it has weak approximation over k} if it is so for any finite $S$. Let
$$\A(S,G) = \prod_{v \in S}G(k_v)/Cl_S(G(k)),~ \A(G) = \prod_v 
G(k_v)/Cl(G(k)),$$ the obstruction (or defect) of weak approximation in $S$ and over $k$,
 respectively, where $Cl$ denotes the closure in the product topology. 
Let $G(k)/R$ (resp. $G(k)/Br$)
denote the group of R-equivalence 
(resp. Br-equivalence) classes of $G$ over $k$. Let
$(\cdot)~^{ \tilde {}}$ be the group 
$\Hom(\cdot, {\bf Q}/{\bf Z}),$ and  $(\cdot)\hat{}$ be the group 
$\Hom(\cdot,{\bf G}_m)$.
$\III(\cdot)$ denotes the Tate - Shafarevich group of $(\cdot)$. We denote also
by $\Br_a~X := \Br_1X /\Br_0 X,$ the {\it arithmetic Brauer group of X}.
By $\H^i(k,G)$ we denote the Galois cohomology of $G$.
\section { Recall of some basic facts from Brauer theory of algebraic groups
[CTS1], [S]}
Let $T$ be a torus defined over field $k$ of characteristic 0 
which is split over a Galois extension 
$K/k$ of group \g. A $k$-torus $N$ is called $induced$, if its
character module $\hat N$ has a {\bf Z}-basis, over which $Gal(\bar k/k)$
acts by permutations. A $k$-torus $S$ is called \g-$flasque$ torus
over $k$ if $\H^{-1}($\h,$\hat S) =0$ for all subgroup
\h~  $\subset$ \g. It is well-known ([CTS1], [V1]) that any
torus $T$ above has a \g -flasque resolution, i.e., an extension
$$1 \to S \to N \to T \to 1 \leqno{(1)}$$
of $T$, where $S$ is a \g -flasque 
$k$-torus, and $N$ is an induced $k$-torus. 
Denote by
$\Br(k,K)$ the kernel of $\H^2(k,{\bf G}_m) \to \H^2(K,{\bf G}_m)$.
The exact sequence (1) induces a homomorphism
$$\H^1(k, \hat S) \stackrel{}{\to} \H^2(k,\hat T), \leqno{(2)}$$
which is injective, since $N$ has trivial 1-cohomology.
\\
One has a cup-product 
$$T(k) \times \H^2(K/k, \hat T) \stackrel{\cup}{\to} \Br(k,K),$$
which defines, via (2), a pairing
$$T(k) \times \H^1(K/k, \hat S) \stackrel{\cup}{\to} \Br(k,K).$$
We have
\\
\\
{\bf 1.1. Theorem.} ([CTS1], Prop. 17 and Corol.)
\\
1) {\it The map $\beta$ defines the
Brauer equivalence  relation over $T(k)$, hence also a map $$\gamma : T(k)/Br \to 
\Hom (\H^1(K/k,\hat S), \Br(k,K)).$$}
\\
2) {\it We have the following anti-commutative diagram}
\\
\\
\hspace*{4cm}$T(k)/R \hspace{0.5cm} \stackrel{p}{\to} \hspace{0.5cm} T(k)/Br$
\\
\\
$\hspace*{4.5cm}\Big \downarrow \delta \hspace{2.7cm} \Big \downarrow \gamma$
\\
\\
\hspace*{3.5cm}$\H^1(k,S)~ \stackrel{\omega}{\to}~ \Hom(\H^1(K/k,\hat S), 
\Br(k,K))$
\\
\\
{\it Here $\delta$ is an isomorphism [CTS1, Theorem 2],
and $\omega$ comes from the cup-product}
$$\H^1(k,S) \times \H^1(k,\hat S) \stackrel{\cup}{\to} \H^2(k,{\bf G}_m).$$
\\
3) $T(k)/Br \simeq \Im(\omega)$ {\it and $T(k)/Br$ is a birational invariant in the class
of k-tori, stably equivalent to T.}
\\
4) {\it If k is a \p -adic local field, the Brauer equivalence on $T(k)$ coincides with
R-equivalence on $T(k)$ and }
$$ T(k)/Br \simeq \H^1(k,\hat S)~^{\tilde{}}.$$
5) {\it If k is a number field, and $\mu$ is the composition map}
$$\H^1(k,\hat S) \stackrel{\Delta}{\to} \prod_v\H^1(k,\hat S) \stackrel{\lambda}{\to}
\prod_v \H^1(k_v,\hat S)$$
{\it then} $T(k)/Br \simeq [\Im(\lambda)/\Im(\mu)]~^{\tilde {}}.$   
\\
\\
{\bf 1.2. Theorem.} ([CTS1], Prop. 19) {\it With above notation, let
 k be a number field. We have the following exact sequences } :
 $$ 0 \to \III(S) \to T(k)/R \stackrel{\rho_T}{\to}
 \prod_v T(k_v)/R \to \A(T) \to 0, \leqno{(R)}$$
$$ 0 \to \A(T) \to \H^1(k,\hat S)~^{\tilde{ }} \to \III(T) \to 0. \leqno{(V)}$$
The following result gives us the group structure on $G(k)/Br$, induced from
that of $G(k)$. Denote by $\Br_e G$ the kernel of the homomorphism
$\Br_1 G \to \Br~k$, defined by specializing at the unit element $e \in G(k)$,
which is isomorphic to $\Br_a G$ ([S], Lemme 6.9). 
\\
\\
{\bf 1.3. Proposition.} ([CTS1], p. 216, [S], Lem. 6.9(1))
{\it Let K be a field 
and G a connected linear algebraic group over K, 
assumed to be reductive if K is
not perfect. Then the pairing}
$$ G(K) \times \Br_eG \to \Br~K$$
{\it is biadditive. In particular, $G(K)/Br$ has a natural group structure
induced from $G(K)$.}   
\\
\\
The following well-known fact (which is a direct consequence
of the Hasse principle for Brauer group of global fields)
was mentioned in [MT] :\\
\\
{\bf 1.4. Proposition.} [MT]
 {\it Let X be a smooth variety defined over a number field k.
Then the restriction map}
$$ X(k)/Br \to \prod_v X(k_v)/Br$$
{\it is injective.}
\\
\\
{\bf 1.5. Remarks.} 1)
 Notice that in Theorem 1.2, we have identified $\III(S)$ with
a subgroup of $T(k)/R$ via the isomorphism $\delta$ of Theorem 1.1, 2).
 The exact sequence (V), which is due to Voskresenski\v i
(see e.g. [V1], [S]), has
been extended to the case of arbitrary connected linear algebraic groups over
number fields by Sansuc [S].
\\
2) We are interested in Brauer equivalence relation for connected
 linear algebraic 
groups, which are rational over algebraic closure $\bar k$ of $k$, hence the 
Br-equivalence and Br$_1$-equivalence are the same.
\\
Our objective is to study
the analogs of the exact sequence (R) in the case of Brauer
and R-equivalence over local and global fields,
for tori in particular, and for connected linear algebraic groups 
in general.
\section {A Brauer relative of exact sequence (R) for algebraic tori}
Let $S$ be a finite set of valuations of a number field $k$, $T$ a $k$-torus,
$T_S := \prod_{v  \in S} T(k_v)$. Denote by $RT(L)$ (resp. $BT(L)$) the set of
elements of $T(L)$ which are R- (resp. Br-) equivalent to 1 in $T(L)$, where
$L$ is a field extension of $k$. Let $RT_S = \prod_{v \in S} RT(k_v)$,
$BT_S = \prod_{v \in S} BT(k_v)$. The following result
 was mentioned in [V2] (which is valid
also for any field $k$ with non-trivial $v$-adic valuations).
\\
\\
{\bf 2.1. Proposition.} $RT_S \subset Cl_S(T(k))$ {\it and is an open subgroup
in $T_S$.}
\\
\\
From above one derives the following
\\
\\
{\bf 2.2. Corollary.} $$A(S,T) \simeq \Coker (T(k)/R \to \prod_{v \in S}
 T(k_v)/R),$$
$$
A(T) \simeq \Coker (T(k)/R \to \prod_v T(k_v)/R).$$
\noindent
{\bf 2.3. Corollary.} $BT_S \subset Cl(T(k)).$
\\
\\
{\it Proof.} If $v$ is a non-archimedean, then Theorem 1.1, (4) tells us that
$BT(k_v) = RT(k_v)$. If $v$ is archimedean, then it is well-known that $T$ is
rational over $k_v$, hence has trivial groups $T(k_v)/Br$ and $T(k_v)/R$,
i.e., $BT(k_v) = RT(k_v)$. \halm
\\
\\
In what follows we identify $T(k)$ with a subgroup of $T_S$ via diagonal 
embedding.
\\
\\
{\bf 2.4. Proposition.} {\it We have}
\\
\begin{tabbing}
1)  $\A(S,T)$ \=$\simeq \Coker(T(k)/Br \to \prod_{v \in S}
T(k_v)/Br).$ 
\\
\\
2) $\A(T) $\>$\simeq \Coker(T(k)/Br \to \prod_v T(k_v)/Br).$
\end{tabbing}
\vspace*{0.6cm}
{\it Proof.}
Notice that
\begin{tabbing}
\hspace*{2cm}$\Coker (T(k)/Br \to T_S/BT_S)$ \=$= T_S/T(k)BT_S$
\\
\\
\>$= T_S/Cl_S(T(k)),$
\end{tabbing}
since $BT_S$ contains $RT_S$ so is also an open subgroup of $T_S$.
So 1) and 2) follow by noticing that for almost all $v$
$$T(k_v)/Br = 1,$$
by [CTS1], p. 205. \halm
\\
\\
We have the following close analog of an exact sequence in [CTS1] (see Theorem
1.2 above) in the case   
of R-equivalence  of algebraic tori over number fields.
\\
\\
{\bf 2.5. Proposition.} $1)$
{\it With above notation we have the following exact
sequence }
$$1 \to T(k)/Br \stackrel {\gamma_T}{\to} \prod_v T(k_v)/Br \to \A(T) \to 1.$$
$2)$ {\it With notation of Theorem 1.1, if
the restriction map $ \H^1(k,\hat S) \to \H^1(k_v,\hat S)$ is surjective
for all v, then  $\Ker (\omega) = \III(S)$.}
\\
\\
{\it Proof.} $1)$ follows directly from Propositions 1.4 and 2.4. 
\\
\\
2) First we show that under the given assumption,
$$\Ker (\omega) \subset \III(S).$$  
Consider the following diagram
\noindent
\\
\\
\hspace*{3cm}$\H^1(k,S) \to \Hom (\H^1(k,\hat S), \Br(k))$
\\
\\
\hspace*{4.1cm}$\Big \downarrow q$ %\hspace{2cm} \Big \downarrow \zeta $
\\
\\
\hspace*{2.4cm}$\prod_v \H^1(k_v,S) \stackrel{\omega'}{\simeq}
 \prod_v \Hom(\H^1(k_v, \hat S), 
\Br(k_v))$
\\
\\
where $\omega'$ is an isomorphism 
by Tate - Nakayama duality. We have, e.g. for $v \not \in \infty$
\begin{tabbing}\\
\hspace*{2cm}$\H^1(k_v,S)$ \=$\simeq \H^1(k_v,\hat S)~\tilde{}$
\\
\\
\>$\simeq \Hom(\H^1(k_v,\hat S), {\bf Q}/{\bf Z})$
\\
\\
\>$= \Hom(\H^1(k_v, \hat S), \Br(k_v))$
\end{tabbing}
Therefore it suffices to show that $\omega' (q(\Ker (\omega))) =0.$
Denote $res_v : \H^1(k,\cdot) \to \H^1(k_v,\cdot)$ the restriction map
of cohomology when passing to completion $k_v$.
Since for $x \in \H^1(k,S)$, $\omega(x)$ is the map
$$ y \mapsto (x) \cup (y), y \in \H^1(k,\hat S),$$
the $v$-component of $\omega'(q(x))$
is given by $$\omega'(q(x))_v : y_v \mapsto (res_v(x)) \cup (y_v),
y_v \in \H^1(k_v, \hat S)$$
 Since the cup-product is compatible with restriction maps,
and each $y_v \in \H^1(k_v, \hat S)$ has the form 
$res_v(y), y \in \H^1(k,\hat S)$,
$$(res_v(x)) \cup (y_v) = (res_v(x)) \cup (res_v(y)) = res_v(x \cup y).$$
Therefore, if $x \in \Ker (\omega)$, then $res_v(x\cup y) =0$ for all $v$,
hence 
 $$\Ker (\omega) \subset \Ker (q) = \III(S).$$
Next we show that $\III(S) \subset \Ker (\omega)$. This is true in 
general without the condition on $res_v$ above. 
We prove that there exists an exact sequence as follows
$$ 0 \to \Ker (\omega)
\stackrel {\alpha}{\to}\III(S)\stackrel {\beta}{\to} T(k)/Br $$
We define the map $\beta : \III(S) \stackrel{i}{\to} T(k)/R \stackrel{p}{\to}
 T(k)/Br$ to be the composite map, where $i$ is the restriction
of the isomorphism
$\delta^{-1}: \H^1(k,S) \to T(k)/R$ (see Theorem 1.1)
to $\III(S)$, and $p$ is the projection.
We show that $\beta$
is a trivial homomorphism, and $\alpha$ is just the identity map.
\\
\\
$a)$ $\Ker(\beta) \subset \Ker(\omega)$. 
We have
$$\Ker (\omega) = \{x \in \H^1(k,S) : (x) \cup (y) = 0, \forall y \in \H^1(k,
\hat S)\}.$$
We have an isomorphism 
$\delta : T(k)/R \simeq \H^1(k,S)$ ([CTS1, Th\'eor\`eme 2]), so  
\begin{tabbing}
\hspace*{2cm} $\Ker(\beta)$ \=
$=\{ x \in \III(S) \subset \H^1(k,S) : \delta^{-1}(x) \in BT/RT \}$
\\
\\
\>$=\III(S) \cap \delta(BT/RT)$
\\
\\
\>$=\{ x \in \III(S) : (x) \cup (y) = 0, \forall y \in \H^1(k,\hat S)\}$
\end{tabbing}
hence $\Ker (\beta) \subset \Ker (\omega)$.
\\
\\
$b)$ $\Im (\beta) = \Ker(\gamma_T) (=0)$.
 Cosider the following commutative diagram
\noindent
\\
\\
\hspace*{2.5cm}$\III(S) \hspace{0.3cm}\stackrel {i}{\to} \hspace{0.5cm}
T(k)/R \stackrel {\rho_T}{\to} \prod_v T(k_v)/R $
\\
\\
\hspace*{3.7cm}$\beta  \searrow \hspace{1cm}  \Big 
\downarrow p' \hspace{1.3cm} \simeq \Big \downarrow q'
\hspace{0.5cm}  $\\
\\
\hspace*{4.7cm} $T(k)/Br \stackrel{\gamma_T}{\to} \prod_v T(k_v)/Br $
\\
\\
If $x \in \Im(\beta)$, then $x = p'(i(s))$, $ s \in \III(S)$. Since $
\rho_T \circ i =0$ by Theorem 1.2,
we have 
\\
\begin{tabbing}
\hspace*{3cm}$q'(\rho_T(i(s)))$ \=$= \gamma_T (p'(i(s)))$
\\
\\
\>$= \gamma_T(x)$
\\
\\
\>$= 0,$
\end{tabbing}
i.e. $ x\in \Ker (\gamma_T)$.
\\
Conversely, if $x \in \Ker(\gamma_T)$, $x = p'(t)$, since $p'$ is surjective.
Then $0 = \gamma_T(p'(t)) = q'(\rho_T(t))$, so
$\rho_T(t) = 0$, since $q'$ is an isomorphism (see Theorem 1.1 (4)). Hence
$t \in \Ker(\rho_T) = \Im(i)$ since the upper row is exact by Theorem 1.2.
\\
Since $\gamma_T$ is injective (see Proposition 1.4), $\beta$ is trivial and
$\alpha$ is an isomorphism. Hence $2)$ is proved. \halm
\\
\\
As a consequence of the above proposition, we have the following
\\
\\
{\bf 2.6. Proposition.} {\it We have the following exact sequence 
connecting the two groups of rational equivalence classes}
$$0 \to \III(S) \to T(k)/R \to T(k)/Br \to 0.$$
{\it  In 
particular, the order of $T(k)/Br$ is equal to $n_T$.}
\\
\\
\Pr
We have (with above notation) the following commutative diagram.
\noindent
\\
\\
\hspace*{2cm}$1 \to \III(S) \to \hspace{0.2cm}T(k)/R \stackrel{\rho_T}{\to}
\hspace{0.15cm} \prod_v T(k_v)/R \to \A(T) \to 1$
\\
\\
$\hspace*{3cm}\Big \downarrow \lambda_T \hspace{1.5cm}\Big \downarrow 
\lambda_T' \hspace{1.5cm} 
\Big \downarrow \simeq \hspace{1.8cm} \Big \downarrow =$
\\
\\
\hspace*{2.8cm}
$1 \hspace{0.5cm}\to \hspace{0.3cm} T(k)/Br \stackrel{\gamma_T}{\to}
 \prod_v T(k_v)/Br \to \A(T) \to 1.$ 
\\
\\
In this diagram, $\lambda_T$ is induced from $\lambda_T'$ and is just the
quotient map. Indeed, we have the vertical
isomorphism "$\simeq$" due to Theorem 1.1, 4), and it is clear that
$$\lambda_T' (\Ker (\rho_T)) \subset \Ker (\gamma_T).$$
Therefore it follows that
$$ T(k)/R / \Ker (\rho_T) \simeq T(k)/Br$$
and we are done.  \halm
\\
\section { Some reductive analogs}
In this section we prove some analogs of results in Section 2 for the case of connected
reductive groups $G$ over number fields $k$. First we recall the following
\\
\\
{\bf 3.1. Proposition.} [T3]
 {\it Let G be a connected linear algebraic groups defined
over a number field k, S a finite set of valuations of k.
 For each $ v\in S$ denote
by $RG_v$ the subgroup of $G(k_v)$ consisting of elements R-equivalent to 1, and
by $RG_S$ the direct product of $RG_v$ for $v \in S$.  
Then $RG_S \subset Cl_S(G(k)).$} 
\\
\\
{\bf 3.2. Proposition.} [T3] {\it Let G, k, S be as above. Then we 
have the following
canonical isomorphisms}
\begin{tabbing}
$1)$ $\A(S,G)$ \=$\simeq~ \Coker (G(k)/R \to \prod_{v \in S} G(k_v)/R).$
\\
\\
$2)$ $\A(G)$ \>$\simeq~ \Coker (G(k)/R \to \prod_v G(k_v)/R).$
\end{tabbing}
 We have the following analog in the case of Brauer equivalence relation.
\\
\\
{\bf 3.3. Theorem.} {\it Let G, k, S be as above. Let $BG_v$ be the subgroup
of $G(k_v)$ consisting of elements which are Br-equivalent to $1$, and
$BG_S$ be the direct product of $BG_v$. Then}
$$BG_S \subset Cl_S(G(k)).$$
{\it Proof.} We follow the proof given in [T3]. We know by [CTS1]
that for a torus $T$ over $k$, $T(k_v)/R = T(k_v)/Br$, hence
it follows that the Theorem holds for tori. Now we may assume that $G$ is not a
torus. Further we just follow the proof given in [T3], where $R$ is replaced by
$Br$ everywhere. \halm
\\
\\
{\bf 3.4. Theorem.}
{\it With notation as above 
we have 
the following
canonical isomorphisms}
\begin{tabbing}
$1)$ $\A(S,G)$  \=$\simeq \Coker (G(k)/Br \to \prod_{v \in S} G(k_v)/Br).$
\\
\\
$2)$ $\A(G)$ \>$\simeq \Coker (G(k)/Br \to \prod_v G(k_v)/Br)$.
\end{tabbing}
{\it Proof.}  
The same as in 2.4, by making use of Proposition 3.3.
\halm
\\
\\
 We need the following technical result.
\\
\\
{\bf 3.5. Proposition.} [O] {\it Let G be a connected reductive group 
defined over a field K. There exists a connected reductive K-group H with
simply connected semisimple part and an induced K-torus Z such that the 
following sequence is exact.}
$$ 1 \to Z \to H \to G \to 1.$$
(Such $H$ is called in the literature also a $z$-$extension$ of $G$ over $K$.)
\\
\\
The relation between the groups of Brauer equivalence classes of $G$ and $H$
is shown in the following statement, where we restrict ourselves only to
the case of a field $k$ of characteristic 0.\\
\\
{\bf 3.6. Proposition.} {\it Let k be a field of
characteristic $0$.}
\\
$1)$ {\it If H is a z-extension of a connected 
reductive k-group G then there is a canonical isomorphism}
$$ H(k)/Br \simeq G(k)/Br.$$
$2)$ {\it If for any torus over k, $T(k)/Br$ is finite, then
 for any connected linear algebraic group G, the group of Brauer
equivalence classes $G(k)/Br$ is finite.
In particular, it is so if
k is  of finite type over {\bf Q} or a local
field of characteristic 0.}
\\
\\
First
we need the following (perhaps well-known for experts, but I do not know of
any reference).
\\
\\
{\bf 3.6.1. Lemma.} {\it Let $X \stackrel{\pi}{\to} Y$ be a morphism of
smooth varieties all 
defined over a field k of characteristic 0, $\pi^* : \Br_1 Y \to 
\Br_1 X $ is an induced
homomorphism. The following diagram is commutative}
\\
\\
\noindent
\hspace*{5cm}$X(k) \times \Br_1 X $
\\
\hspace*{7.5cm} $\searrow$
\\
\hspace*{5.5cm}$\Big \downarrow \pi \hspace{0.6cm}
\Big  \uparrow \pi^* \hspace{0.7cm}
 \Br~ k$
\\
\hspace*{7.5cm} $\nearrow$ 
\\
\hspace*{5cm}$Y(k) \times \Br_1 Y $
\\
{\it i.e. the pairing
$$X(k) \times \Br_1 X \to \Br~ k$$ is functorial in $X$, i.e., given
$x \in X(k)$, $b \in \Br'X$, where $\Br'$ denotes the usual Brauer group of
$X$, which we identify with the cohomological one (see [S, Sec. 6]), then
}
$$\pi^* (b)(x) = b(\pi (x)).$$
\Pr  Let $x \in X(k)$, $y = \pi(x) \in Y(k)$. Denote by ${\cal O}_{\bar X,x},
{\cal O}_{\bar Y, y}$ the local ring of $\bar X$ (resp. $\bar Y$) at $x$ (resp. $y$). One has the following commutative diagram
\\
% (x,y) = coordinates ! The bigger y the higher the point \mbox{...}
\begin{center}
\begin{picture}(100,70)(0,45)
\put(-20,90) {${\cal O}_{\bar Y,y}$}
\put(21,94){\vector(1,0){50}}
\put(51,97){$\pi^*$}
\put(72,90){ \mbox{${\cal O}_{\bar X,x} $}}
\put(75,42){ \mbox{$\bar k $}}
\put(88,85){\vector(0,-1){29}}
\put(90,69){$\pi'$}
\put(12,87){\vector(2,-1){65}}
\put(40,55){$\pi''$}
\end{picture}
\end{center}
We denote
\begin{center}
$\gm~ := Gal(\bar k/k),$
$ (f_{s,t}) \in Z^2(\gm , {\cal O}_{\bar Y,y}^*)$
\end{center}
the absolute Galois group of $k$
and a 2-cocycle representative of $b \in \Br_1 Y$, respectively.
Then by [S], Lemme 6.2,
$b(y)$ is the class $[f_{s,t}(y)]$ in $\Br~ k$. This 2-cocycle gives rise
to a 2-cocycle $[f_{s,t} \circ \pi] \in Z^2(\gm , {\cal O}_{\bar X, x}^*)$
which is nothing else than a representative of
$a = \pi^* (b)$. Then $a(x)$ is just the class of
$$[(f_{s,t} \circ \pi)(x)] = [f_{s,t}(\pi (x))] = [f_{s,t}(y)] \in \Br~ k.$$
Therefore $$\pi^* (b)(x) = b(\pi (x))$$ as required. 
The proof of the lemma is complete.  \halm
\\
\\
{\it Proof of Proposition 3.6.}
\\
\\
$1)$ Let $\pi : H \to G$ be the projection. It induces an epimorphism 
$$H(k) \to G(k)$$
hence also an epimorphism
$$\pi' : H(k)/Br \to G(k)/Br.$$
We show that $\pi'$ is injective.
Since $Z$ is an induced $k$-torus, it is well-known that there
is a $k$-section 
$$ i : G \to H, \pi \circ i = id_G.$$
By Lemma 3.6.1, we have the following commutative
diagram
\\
\\
\hspace*{5cm}$H(k)\hspace{0.5cm} \times \hspace{0.5cm} \Br_1 H $
\\
\hspace*{8.5cm} $\searrow$
\\
\hspace*{5cm}$i \Big \uparrow \Big \downarrow \pi \hspace{1.2cm}  \pi^*
\Big \uparrow 
\Big \downarrow i^* \hspace{0.8cm} \Br~ k$
\\
\hspace*{8.5cm} $\nearrow$ 
\\
\hspace*{5cm}$G(k) \hspace{0.5cm}\times \hspace{0.5cm} \Br_1 G $
\\
\\
Since $\pi \circ i =id_G$, it follows that we have
$$(\pi \circ i)^* = i^*\circ \pi^* = id_{Br_1G}.$$
Therefore $i^*$ is {\it surjective.} Now assume that
$h \in H(k)$ such that
$$\pi(h) \cup \hat g = 0, \forall \hat g \in \Br_1G.$$
Let $g = \pi(h)$. Then by Lemma 3.6.1 we have
$$ i(g) \cup \hat h =  \pi(i(g)) \cup i^*(\hat h) = 0, \forall
\hat h \in \Br_1H$$
since $i^*$ is surjective. This implies that
$i(g) \in BH(k)$. Since $$\pi(h) = g = \pi (i(g)),$$
one deduces that $ h \equiv i(g) (mod. Z(k))$.
Since $Z$ is $k$-rational, it follows that
$Z(k) \subset BH(k)$, hence we are done.
\\
(In the case 2) we can prove our statement
by using 2). 
 Indeed, we have a birational equivalence
$$H \simeq G \times Z, $$
and this induces a bijection (see [CTS1], Section 7)
$$H(k)/Br \simeq (G(k)/Br \times Z(k)/Br) \simeq G(k)/Br.$$
since $H(k)/Br$ is finite by $2)$, $\pi'$ is injective, hence
also an isomorphism.)
\\
\\
$2)$ One reduces easily to the case where $G$ is a connected reductive group.
Take a $z$-extension $H$ of $G$ as above. We will show that
$H(k)/Br$ is finite. Let $H = S\tilde G$, where $\tilde G$ is
a  simply connected semisimple group and $S$ a central torus of $H$. We have
the following exact sequence of $k$-groups
$$ 1 \to \tilde G \to H \to T \to 1$$
where $T$ is a torus. One derives from [S], Corollaire 6.11, the following
exact sequence
$$ \Pic T \to \Pic H \to \Pic \tilde G \to \Br_a T \to 
\Br_a H \to \Br_a \tilde G.  \leqno{(3)}$$
By [S], Lemme 6.9 (iv), we know that $\Pic(\tilde G)$ and $\Br_a \tilde G$
are trivial, so we obtain from (3) the following isomorphism
$$\Br_a T \simeq \Br_a H.$$
Since $\Br_1 T \hookrightarrow \Br_1 H$ ([S], Corollaire 6.11),
 it follows that $$\Br_1 T = \Br_1 H. \leqno{(4)}$$ 
(One can use also Prop. 6.10 (loc.cit) to get this equality.)
\\
\\
We continue the proof of the proposition. By Lemma 3.6.1 we have the following 
commutative diagram
\\
\\
\noindent
\hspace*{5cm}$H(k) \times \Br_1 H $
\\
\hspace*{7.5cm} $\searrow$
\\
\hspace*{5.5cm}$\Big \downarrow \pi \hspace{0.6cm} \Big\uparrow \pi^* \hspace{0.7cm} \Br~ k$\\
\hspace*{7.5cm} $\nearrow$ 
\\
\hspace*{5cm}$T(k) \times \Br_1 T $
\\
\\
Assume that $h \in H(k)$ with $\pi(h) \in BT(k)$. Then we have
$$\pi^*(a)(h) = a(\pi (h)) = 0$$
for all $a \in \Br_1 T,$ and from (4) we derive that 
$h \in BH(k)$. Thus we have proved that $\pi$ induces an injective
homomorphism 
$$H(k)/Br \to T(k)/Br. $$
Therefore, the first statement of 2)
follows. Now if $k$ is a field of finite type over {\bf Q} or a local
field, then by [CTS1], Corollaire 1, p. 217, $T(k)/Br$ is finite. Therefore
$H(k)/Br$, and {\it a fortiori} $G(k)/Br$, is also finite.
\halm
\\
\\
{\bf 3.6.2. Remark.} It was proved in [G2] that if
$k$ is a number field then $G(k)/R$ is finite, hence $G(k)/Br$
is also. Here we did not use the result of Gille in the proof.
\\
\\ 
{\bf 3.6.3. Corollary.} (of the proof.) {\it Let}
$$ 1 \to A \to B \stackrel{\pi}{\to} C \to 1$$
{\it be an an exact sequence of connected linear algebraic groups 
defined over a 
field k of characteristic 0, where A has trivial Picard and arithmetic Brauer
groups. Then the canonical homomorphism $B(k)/Br \to C(k)/Br$ is injective.
In particular, if k is a non-archimedean local field and A is simply connected
then we have}
$$ B(k)/Br \simeq C(k)/Br.$$
\Pr  Only the second part requires a proof, which follows from Kneser's
Theorem on the triviality of $\H^1$ of simply connected groups. Indeed, from the exact
sequence of cohomology we see that $\pi$ is surjective on $k$-points, thus
gives a surjective map
$$B(k)/Br \to C(k)/Br,$$
which is also bijective.
\halm
\\
\\
Before we give the formulation of one of main results,
 we recall some definition 
and notation of Section 2.4.
\\
For a reductive group $H$, we call {\it torus quotient} of $H$ the factor
group of $H$ by its semisimple part $[H,H]$. Any connected linear 
algebraic group $G$ over a field $k$ of characteristic 0 is a semidirect
product $G = L R_u(G)$, where $L$ is a Levi (reductive) subgroup of $G$,
and $R_u(G)$ is the unipotent radical of $G$. $L$ is unique up to conjugacy
by elements from $G(k)$ and by convention, we call the 
{\it semisimple part of G} the semisimple part of some fixed Levi
$k$-subgroup of $G$. 
Given a torus $T$ defined over a number field $k$ we denote by $V(T)$ a smooth
compactification of $T$ over $k$ and by $S$ the Neron - Severi $k$-torus of $T$,
which is by definition the
dual to the Picard group of $\bar V(T) (= V(T) \times \bar k)$,
 $\hat S \simeq \Pic(\bar V(T))$. The first Galois
cohomology of $S$, which depends on $T$,
does not depend on the chosen smooth compactification, and so is the
Shafarevich - Tate group $\III(S)$.
By Proposition 2.5 we have the following exact sequence
$$1 \to T(k)/Br \stackrel{\gamma_T}{\to} 
\prod_v T(k_v)/Br \to \A(T) \to 1.$$
A generalization  of this sequence is given by the following theorem.
\\
\\
{\bf 3.7. Theorem.} 
{\it Let G be a connected linear algebraic group defined over a number field k. Then
we have  the following commutative diagram, all rows of
which are exact sequences}
\\\
\\
\noindent
\hspace*{3.3cm}$1 \to G(k)/Br \stackrel{\gamma_G}{\to}\prod_v G(k_v)/Br
\to \A(G) \to 1$
\\
\\(5)
\hspace*{4cm}$\Big \downarrow p \hspace{2cm}
\Big \downarrow q \hspace{2cm}\Big \downarrow r$
\hspace{2.7cm}
\\
\\
\hspace*{3.1cm}
$1 \to T(k)/Br \stackrel{\gamma_T}{\to} \prod_v T(k_v)/Br \to \A(T) \to 1$
\\
\\
{\it where T is the torus quotient of any z-extension H of the reductive
part of G, and all vertical maps are (functorial) isomorphisms
(including local components ones).
 In particular,
the image of $G(k)/Br$ via $\gamma_G$ is a finite group of order 
$n_T$.} 
\\
\\
\Pr  
The exactness of the above sequences follows from Propositions 2.5, 1) and 3.4.
By Proposition 3.6 (1), there is a canonical (functorial) isomoprhism
$G(K)/Br \simeq H(K)/Br$ for any extension field $K/k$. Therefore it suffices
to prove theorem 3.7 for $H$.
\\
We have the following commutative diagram (see the notation above).
\noindent
\\
\\
\hspace*{3.3cm}$1 \to H(k)/Br \stackrel{\gamma_H}{\to}\prod_v H(k_v)/Br
\to \A(H) \to 1$
\\
\\
\hspace*{4.5cm}$\Big \downarrow p \hspace{2cm}\Big
\downarrow q \hspace{2.4cm}\Big \downarrow r$
\hspace{2.7cm}
\\
\\
\hspace*{3.1cm}
$1 \to T(k)/Br \stackrel{\gamma_T}{\to}~ \prod_v T(k_v)/Br~ \to~ \A(T) \to 1$
\\
\\
where $p, q, r$ are natural maps, induced from the
projection $pr : H \to T$.
By Corollary 3.6.3,
$p$ is injective and $q$ and all its local components of
are isomorphisms.
\\
Next we need the following
\\
\\
{\bf 3.8. Lemma.} {\it With above notation, we have a canonical isomorphism 
of finite groups}
$$ \A(H) \stackrel{r}{\simeq} \A(T).$$
\Pr  Consider the following commutative diagram
%\newpage
\noindent
\\
\\
(6)\hspace*{3cm}
$\matrix{H(k)& \to& T(k)& \stackrel{\delta}{\to}& \H^1(k,\tilde G)\cr
&&&& \cr
\Big \downarrow \alpha &&  \Big
\downarrow \beta && \Big  \downarrow \gamma \cr
&&&& \cr
H_S &\to& T_S &\stackrel{\delta}{\to}& \prod_{v \in S} \H^1(k_v,\tilde G) \cr}$
\\
\\
\\
where we take $S$ a sufficiently large finite set of valuations of $k$ containing
all the archimedean ones, such that
$$ \A(H) = H_S /Cl_S(H(k)),$$
$$\A(T) = T_S/Cl_S(T(k)).$$
It is a general and well-known fact (see, e.g. [T3] for a discussion with references)
 that for any linear algebraic group $P$ over $k$,
$Cl_S(P(k))$ is an open subgroup of $P_S$.
Also any connected linear algebraic group satisfies weak approximation with respect to 
the set $\infty$ of archimedean valuations.
\\
Let $t_S \in T_S$, $t_S = (t_{\infty}, t_f)$, where $t_{\infty}$
(resp. $t_f$) is the $\infty$-(resp. finite) component of $t_S$. Let
$t_n \in T(k)$ such that $lim_n~t_n = t_\infty.$
Then for $n$ large enough, the element
$(t_n,t_n)$ is very close to $(t_\infty, t_n)$ in the $S$-adic (product) topology,
i.e., $(t_n^{-1}t_\infty,1)$ is approaching 1 in $T_S$. Since 
$Cl_S(T(k))$ is open, there is $N$ such that if $n>N$, then
$$(t_n^{-1}t_\infty, 1) \in Cl_S(T(k)),$$
i.e., $$(t_\infty,t_n) \in Cl_S(T(k)).$$
Since $t_S = (t_\infty,t_n) (1,t_n^{-1}t_f)$, it follows that
\\
\\
(7)~  {\it  each coset of
$T_S/Cl_S(T(k))$ has a representative from $1 \times T_{S-\infty}$.} 
\\
\\
Since $\H^1(k_v,\tilde G)$ is trivial
for $v$ non-archimedean by Kneser's Theorem [Kn2], it follows that
$$T_{S - \infty} = \pi (H_{S-\infty})$$
and from (7) we derive that the natural homomorphism
$$\A(S,H) = \A(H) \stackrel{\pi}{\to} \A(T) = \A(S,T)$$ is surjective.
\\
Next we show that $\pi$ is injective. Let $h_S \in H_S$ such that 
$\pi(h_S) \in Cl_S(T(k))$. Then 
$$\pi(h_S) = lim_n~ t_n , t_n  \in T(k),$$
hence from the commutative diagram (6) we derive
$$1 = \delta \pi(h_S) = lim_n~ \delta(t_n).$$
(Notice that here one endows $\H^1(k,\tilde G)$, $\H^1(k_v,\tilde G)$
with discrete
topologies, and one checks readily that all maps in the 
diagram (6) are continuous.)
Since $\prod_{v \in S} \H^1(k_v,\tilde G)$ is finite, there is $N_1$ such that 
if $n > N_1$ then $\delta(t_n) =1$. Let $h_n \in H(k)$, such that
$\pi(h_n) = t_n$. Then $lim_n~ \pi(h_n^{-1}h_S) = 1$ since
$Cl_S(H(k))$ is open in $H_S$, and we deduce that
$$h_n^{-1} h_S \in \pi(\tilde G_S) Cl_S(H(k))$$
for $n$ large.
\\
Since $\tilde G$ has weak approximation property (Kneser - Harder Theorem,
see  e.g. [S]), we have
$\pi(\tilde G_S) \subset Cl_S(H(k))$, hence $h_S \in Cl_S(H(k))$ as required.
Hence $\A(H) \simeq \A(T)$ and the lemma is proved. \halm
\\
\\
{\it Continuation of the proof of Theorem 3.7.} In the diagram (5) we know that
$r$ is an isomorphism by Lemma
 3.8, $q$ is an isomorphism and $p$ is injective by Corollary 3.6.3.
It follows that $\gamma_H$ and $\gamma_T$ have isomorphic images. Therefore
$$ H(k)/Br \stackrel{p}{\simeq} T(k)/Br.$$
In particular,
the order of $G(k)/Br$ is equal to $n_T$.
\halm
\\
\\
We have the following result which is an analog of [CTS1, Prop. 18]
and [S, Th. 3.3]. Let $G$ be a connected 
linear algebraic group defined over number field $k$, $H$ be a $z$-extension
of a Levi subgroup of $G, T$ be torus quotient of $H$, which is split over
a finite extension $K$ of $k$. Denote by $S$ the Neron - Severi torus of
$T$, $V_0$ the (finite) set
of all valuations of $k$, such that their extensions to $K$
have {\it non-cyclic} decomposition groups. Then for any finite
set $W$ of valuations of $k$ we have the following formulas  
\\
\\
{\bf 3.9. Corollary.} (of the proof of Lemma 3.8) {\it There 
are canonical isomorphisms of finite groups } 
$$\A(W,G) \simeq \A(W,T) \simeq \Coker (\H^1(k,S) \to \prod_{v \in W} 
\H^1(k_v,S),$$
$$\A(G) \simeq \A(T) \simeq \Coker (\H^1(k,S) \to \prod_{W \cup V_0}
\H^1(k_v,S).$$
\Pr
The proof of the first isomorphisms
in these chain of isomorphisms 
 follows directly from the proof of Lemma 3.8 above. The last
ones related with $S$ are deduced from Theorem 1.2 (well-known).
\halm 
\section {Some variations and applications}
In this section we consider some applications of results obtained in previous
sections, and also of those obtained in [T3], where we have not given any.
We keep our notation as above. First we derive from Proposition 3.3 the following.\\
\\
{\bf 4.1. Proposition.} {\it Let S be a finite set of valuations of k and G a connected
linear algebraic group over k. If G has trivial group $G(k_v)/Br$ of Brauer equivalence 
classes  for all $v \in S$, then G has weak approximation property in S.}
\\
\\
{\bf 4.2.} Conversely, assume that $G$ does not have weak approximation property
with respect to some $v$ (necessarily non-archimedean). Then $G$ has non-trivial
group $G(k_v)/Br$ by Proposition  3.3, hence also non-trivial group $G(k_v)/R$. Therefore $G$
is not stably rational {\it over $k_v$}, and a fortiori, over $k$.
\\
\\
{\bf 4.2.1. Remarks.}
 $1)$ Usually the counter-examples to weak approximation over 
number fields $k$
serve as examples of linear algebraic $k$-groups, which
are non stably rational over $k$ only.  
The statement 4.2 above shows that these examples are, in fact, stronger in the sense that
they serve also as examples of non stably rational groups over bigger fields (say $k_v$).
The reader may consult a variety of examples in [S].  
\\
\\
$2)$ We mention the following one
of the main results due to Sansuc [S], Corollaire 9.7,
in the case $G$ has no simple component of type $\E_8$ (and also goes back
to Voskresenski in 
torus case).
Since the Hasse principle is also holds for $\E_8$ by Chernousov, the following
holds.
\\
\\
{\bf 4.2.2. Theorem.} ([S], Cf.Thm 9.5.)
 {\it If G is a connected linear algebraic group, defined
over a number field k, then we have the following exact sequence}
$$ 0 \to \A(G) \to \H^1(k,\Pic \bar V(G))~\tilde{}
\to \III(G) \to 0.$$
{\it In particular, if $\H^1(k, \Pic \bar V(G))$ is trivial,
 then G has weak approximation over k and satisfies
Hasse principle for $\H^1$.}
\\
\\
One may ask, by comparing with 4.2, if we have a similar
situation if we assume that
$\A(G) = 0$, and $\III(G) \ne 0$. 
However it is not true as the following classical example shows.
\\
\\
{\it Example.} Let $a,b \in {\bf Z}$ (the integers), and let $K = {\bf Q}(\sqrt a,
\sqrt b),$ be a bi-quadratic extension of {\bf Q}, 
where {\bf Q} denotes the rational numbers. Denote by $T(a,b) = 
R^{(1)}_{K/{\bf Q}}
({\bf G}_m)$. Then (see [CTS1], Prop. 7, or [V1], p.157) we have
$$\H^1({\bf Q}, \Pic\bar V(T)) = {\bf Z}/2 {\bf Z}.$$
If we choose $a,b$ such that all the decomposition groups for $K$ are cyclic
then it is known (by Serre) that $\A(T) =0$. For example,
$a = 5, b = 29$ satisfy
this condition. However, one checks that $T(5,29)$ is rational over all completions
of {\bf Q}, but $\III(T(a,b)) = {\bf Z}/ 2{\bf Z}.$
\\
\\ 
In the next result we consider some applications to weak approximation
in semisimple groups defined over number fields $k$.
\\
\\
{\bf 4.3. Theorem.} {\it Let G be a semisimple k-group such that G is of inner
type over $k_v$ for all $v \in S$.Then G has weak approximation
 over k with respect
to S. In particular, if G is of inner type over k, it has weak approximation
over k.}
\\
\\
The theorem follows from the following
\\
\\
{\bf 4.4. Proposition.} {\it If as a group over
$k_v$, G is an inner type 
then G has trivial group $G(k_v)/R$.}
\\
\\
First we need the following result due to Gille [G1], Prop. 2.3.
\\
\\
{\bf 4.4.1. Proposition.} [G1] {\it Let}
$$1 \to F \to G_1 \stackrel{\lambda}{\to} G_2 \to 1$$
{\it be an isogeny of connected reductive groups, all defined over a field
k of characteristic $0$ and $C_{\lambda}(k) = G_2(k)/\lambda (G_1(k))$. Then the 
following sequence of groups is exact}
$$G_1(k)/R \to G_2(k)/R \stackrel{\delta_R}{\to} \H^1(k,F)/R ,$$
{\it where the image of $\delta_R$ is the factor group
$C_{\lambda}(k)/R$ and 
the R-equivalence relation on $C_{\lambda}(k)~
( \hookrightarrow  \H^1(k,F))$ is induced from
that on $\H^1(k, F)$.}
\\
\\ 
{\bf Note.} In the case that $\H^1(k,G_1)$ is trivial, we identify 
$C_{\lambda}(k)$
with $\H^1(k,F)$ and also write the above exact sequence in the form
$$G_1(k)/R \to G_2(k)/R \to \H^1(k,F)/R \to 0.$$
{\it Proof of Proposition 4.4.}  We distinguish two cases.
\\
\\
$1)$ $v \in \infty$. It is well-known that any connected linear algebraic 
group $G$ over {\bf R} has trivial group of R-equivalences. Here is a short
indication of proof. One reduces easily to proving that any semisimple
element $s \in G({\bf R})$ is R-equivalent to 1. But this follows from
the fact that $s$ is belong to some torus defined over {\bf R}, and any
such torus is rational over {\bf R}.
\\
\\
$2)$ $v$ {\it is non-archimedean.} Let $G_s$ be a $k_v$-split form of $G$, which
is obtained from $G$ by an inner twist. Denote by $F$ the fundamental group
of $G$, which is the same for $G_s$.  One
can check that simply connected groups have trivial
groups of R-equivalence classes, we
 have by Proposition 4.4.1 ([G1], Prop.2.3) the following exact sequences
$$0 \to G_s(k_v)/R \to \H^1(k_v,F)/R \to 0, \leqno{(8)}$$      
$$0 \to G(k_v)/R \to \H^1(k_v,F)/R \to 0. \leqno{(9)}$$
Here we identify $\H^1(k_v,F)$ with the factor group
 $G_s(k_v)/\pi(\tilde G_s(k_v))$
(resp. $G(k_v)/\pi(\tilde G (k_v))$), where $\pi : \tilde G_s \to G_s$
(resp. $\pi : \tilde G \to G$) is the simply connected covering of $G_s$
(resp. $G$) and take the factor group modulo the rational relation
induced on $\H^1(k_v,F)$.
 We have also used the Kneser Theorem on the vanishing of
$\H^1(k_v,\tilde G_s)$ and $\H^1(k_v,\tilde G)$.
Since $G_s$ is rational over $k_v$, the second group in (8) is trivial,
therefore  by (9) the group  $G(k_v)/R$ is also trivial. \halm
\\
\\
Thus by this proposition, and by Proposition 3.3, $G$ has weak
approximation with respect to $S$ and Theorem 4.3 is proved.  \halm
\\
\\
{\bf 4.5. Remarks.}
1) One cannot simply drop the inner type assumption, since there
are examples of semisimple quasi-split groups over number fields which do
not have weak approximation property. First examples of such groups were
given by Serre (see [Kn1] and [S] for more information).
\\
\\
2) In the case of number field,
this result also extends the previously known (but more general)
result by Harder,
namely we have 
\\
\\
{\bf 4.5.1. Corollary.} (Harder [H1]) {\it If a semisimple group G defined over
a number field k is split over $k_v$ for all $v \in S$, then G has weak
approximation in S.}
\\
\\
Now we apply our results to give new proofs (and also
discuss some extension)
of some results due to
Harder and Sansuc. In the following theorem, the first result is due to
Harder ([H2], Satz 2.2.3) and the second is due  to Sansuc ([S], Cor. 5.4).
The third, in the case that the given group is semisimple 
and split over a metacyclic extension 
of the given number field, is also due to Sansuc.
\\
Let $\pi : \tilde G \to G$ be a central isogeny of semisimple groups,
all defined over a field $k$, where $\tilde G$ is simply connected covering of
$G$.
$\pi$ is called a {\it normal isogeny} (after Harder [H2]) if $\mu :=\Ker \pi$
can be embedded into an induced $k$-torus $M$, such that $M/\mu$
is also an induced $k$-torus. One can show, for example, that adjoint groups   
$G$ have normal isogenies.
\\
\\
{\bf 4.6. Theorem.} {\it The following groups have weak approximation property
over number fields.}
\\
$1)$ ([H2]) {\it Semisimple groups which are images of normal isogenies};
\\
$2)$ ([S]) {\it Absolutely almost simple groups;}
\\
$3)$ {\it Inner forms of connected reductive groups which are split over
a metacyclic extension of $k_v$ for all non-archimedean v. Moreover, 
two connected
reductive groups, which are innner form of each other
 have the same group of R-equivalence classes over local 
non-archimedean fields.} 
\\
\\
\Pr  $1)$
Let $\pi : \tilde G \to G$ be a normal isogeny defined over a number
field $k$, and $S$ any finite set of valuations of $k$. We show that for all
$v \in (S - \infty)$, $G(k_v)/R $ is trivial. Indeed, let $\mu = \Ker \pi$,
$M$ be an induced $k$-torus, such that $M/\mu$ is also an induced 
$k$-torus. As above (see Proposition 4.4.1), we have the following exact sequences
$$ \tilde G(k_v) \to G(k_v) \to \H^1(k_v, \mu) \to 0, $$
$$\tilde G(k_v)/R \to G(k_v) /R \to \H^1(k_v,\mu)/R \to 0. $$
One can show easily that $\tilde G(k_v)/R$ is trivial. (Here is a short argument. 
One reduces to almost simple case. If $G$ is isotropic, then 
it is well-known
that $G(k_v)$ has no nontrivial normal subgroup, i.e. $G(k_v) = RG(k_v)$,
since $RG(K_v)$ is a normal Zariski dense subgroup of $G(k_v)$.
Otherwise, $G$ is of inner type $\A_n$ by a result of Kneser, 
and in this case the result
is well-known.)
\\
Hence we have an isomorphism
$$ G(k_v)/R \simeq \H^1(k_v,\mu)/R.$$
By considering similar exact sequences 
$$ 1 \to \mu \to M \to M' \to 1,$$
$$ M(k_v)/R \to M'(k_v)/R \to \H^1(k_v,\mu)/R \to 0,$$ 
where $M, M'$ are induced tori and using that $M,M'$ are rational, 
we get that $\H^1(k_v,\mu)/R$ is trivial, and so is $G(k_v)/R$. Therefore $G$
has weak approximation over $k$ by Proposition 3.3.
\\
In the case $G$ is an adjoint group,
we can use a direct argument as follows.
 We may use the following (easy to show) fact : adjoint groups over
local \p-adic fields are rational. This was mentioned in the preprints
[T1-2]. One may also argue as follows. Let $G_q$ be a quasi-split inner
form of $G$ defined over $k$, and $F$ be its fundamental group. By using 
the same argument we have in the proof of Proposition 4.4, one concludes
that $G(k_v)/R$ is trivial for all $v$. Therefore by Proposition 3.3,
$G$ has weak approximation over $k$.
\\
\\
$2)$ Let $G$ be an 
(absolutely) almost simple $k$-group. We want to show that
for any finite set $S$ of valuations of $k$, $G(k_v)/R$ is trivial for
all $v \in S$. [In fact, one can show a much stronger result in this case :
see the paper by Chernousov and Platonov : 
The rationality problem
of simple algebraic groups, C. R. Acad. Sci. Paris 322 (1996), 245 - 250, 
 which have many results overlapped with results of
[T2], where also other results were mentioned)
\\
{\it If v is a non-archimedean valuation, then 
G is rational over $k_v$ if
G is not of type $\A_n$.}]
\\
Here we can use the following simple argument as follows. 
Let $G_q$ be an almost simple quasi-split inner form of $G$. As in the case
1) we are reduced to proving the statement for quasi-split groups. Assuming
that $G$ is not split, then $G$ is of type $\A_n$, $\D_n$,
$\E_6$, or $^{3,6}\D_4$. Let $T$ be a maximal $k$-torus of $G$ 
containing a maximal
$k_v$-split torus of $G$.
If $G$ is not of trialitarian type, then we know by Tits [Ti] that $T$ is split
by a quadratic extension of $k_v$. The structure of tori split over a quadratic
extensions are well-known : they are direct product of  groups
of type ${\bf G}_m$, $\R_{K/k_v}({\bf G}_m)$, or $\R_{K/k_v}^{(1)}({\bf G}_m)$
where $K/k_v$ is a quadratic extension of $k_v$. In particular they are rational
over $k_v$,  and so is $G$ by Bruhat decomposition. 
In the trialitarian
case one proves in the same way 
that maximal tori containing a maximal split torus are rational.
Thus by Proposition 3.3,
$G$ has weak approximation in $S$ for any $S$, thus also over $k$.  
\\
\\
3) $a)$
First we show that if
$G$ is a connected reductive group which is split over 
metacyclic extension $l_v$ of $k_v$ for each non-archimedean $v \in S$ then
$G$ has weak approximation with respect to
any such finite set $S$ of valuations. In fact we prove the following 
stronger result.
\\
\\
{\bf 4.7. Proposition.} {\it If G is a connected reductive group defined
over non-archimedean $k_v$ and split over a metacyclic extension $l_v$ of $k_v$ then
$G(k_v)/Br$ is trivial.}
\\
\\
\Pr  
It can be shown that there exists a maximal $k_v$-torus
$T$ of $G$ which is split over $l_v$.
(And in the case of number field $k$, one can show that there exists
a maximal $k$-torus $T$ such that $T$ is $l_v$-split for $v \in S$.)
 Let $H$ be a $z$-extension
of $G$,
$$ 1 \to Z \to H \stackrel{\pi}{\to} G \to 1.$$
Let $T_H$ be the maximal $k_v$-torus of $G$ such that $T_H$ is mapped
onto $T$ via $\pi$. 
Let $\tilde G$ be the simply connected covering of the semisimple part
$G'$ of $G$, and let $\tilde T$ be the maximal $k_v$-torus  
of $\tilde G$ which is mapped into $T$ via the composite map
$$ \tilde G \to G' \to G.$$
We have the following exact sequences of tori.
$$ 1 \to Z \to T_H \to T \to 1 , \leqno{(10)}$$
$$ 1 \to \tilde T \to T_H \to T_0 \to 1.\leqno{(11)}$$
It is clear from (10), (11)
that $\tilde T$ is also split over $l_v$ for all $v \in S$. By
[CTS1], Corollaire 3, p. 200, we have
$$\tilde T(k_v)/R = T(k_v)/R = \{1\}, \forall v \in S,$$
therefore $T_H(k_v)/R = \{1\}$, $\forall v \in S$.
\\
Now we consider a maximal $k_v$-split torus $T_1$ of $G$. Then
$$Z_G(T_1) = T' H',$$
where $T'$ is the connected center and
$H'$ is a semisimple $k_v$-group, anisotropic over $k_v$.
If $H'$ is trivial, i.e., $G$ is quasi-split over $k_v$, then the torus $T'$
is split over  metacyclic extension $l_v$, so has trivial group
of R-equivalence classes by [CTS1], Corollaire 3, p. 200, and so is $G$,
since $G$ and
$Z_G(T_1)$ are birationally equivalent over $k_v$ (using Bruhat decomposition).
Therefore $G(k)v)/Br$ is trivial also.
One may therefore assume that $H'$ is non-trivial, and
by replacing $G$ by $Z_G(T_1)$, one may
assume that $G$ has semisimple part $G'$ anisotropic over
$k_v$.
\\
By using a consequence of the Kneser's Vanishing Theorem [Kn2],
we see that $G'$ is necessarily a product of almost
simple $k_v$-factors of type $^1\A$, which may be taken to be
absolutely almost simple. So we have
$$\tilde G = H_1 \times \cdots \times H_t,$$
where $H_i$ is simply connected of type $^1\A_{n_i -1}$ for all $i$.
\\
We now recall the construction of the $z$-extension $H$ of $G$. Let
$G = G' P$, where $P$ is a $k_v$-torus. Let $F = G' \cap P$,
$F_1 = \{(f,f) : f \in F\}$ and we have
the following exact sequence
$$ 1 \to  F_1 \to G' \times P \to G \to 1. $$
By taking the composite of two isogenies 
$$ \tilde G \times P \to G' \times P \to G,$$
we have an isogeny
$$1 \to F' \to \tilde G \times P \to G \to 1.$$ 
Thus one sees that since $\tilde G$ is of inner type (in fact the product of
groups $\SL$), the group $F'$ can be embedded into a split torus
$Z$ defined over $k_v$. Then we take $$F_1' : =
\{(f,f) : f \in F'\}$$
and take
$$H = (Z \times (\tilde G \times P))/F_1'.$$
From the very construction it follows from (10), (11)
that $T_H$ is split over $l_v$. Since $\tilde T$ and $T_H$ are split over
$l_v$, the same holds for $T_0$. Therefore $T_0(k_v)/R$ is trivial
by [CTS1], p. 200. Since $T_0(k_v)/R = T_0(k_v)/Br$ ([CTS1], p. 217)
and $H(k_v)/Br = T_0(k_v)/Br$ by Corollary 3.6.3 and the proof given there,
we see that $H(k_v)/Br = G(k_v)/Br = \{1\}$.
The proof of 4.7 is complete. \halm
\\
Now we see that $G(k_v)/Br =1$ 
 for all $v \in S$.
Thus $G$ has weak approximation for any given finite set $S$, which means
that $G$ has weak approximation over $k$.
\\
\\
3) $b)$ Now we assume that $G_1$ is an inner form of a group $G$.
 First we prove the last statement in 3) of the theorem.
\\
We need the following very useful
\\
\\
{\bf 4.8. Lemma.} (Ono - Sansuc [S]) {\it Let G be a connected reductive group 
defined over a field k. There exists a number n, induced 
k-tori T and T' such that
we have the following central k-isogeny}
$$ 1 \to F \to \tilde G^n \times T' \to G^n \times T \to 1,$$
{\it where $\tilde G$ is the simply connected covering of the semisimple part
G' of G.}
\\
\\
The finite covering of an algebraic  group
by a direct product of simply connected group with an induced
torus (such as $\tilde G^n \times T' \to G^n \times T$ above)
is called after Sansuc  a
{\it special covering}.
It is obvious that to prove our statement we may assume that the group $G_1$
itself has a special covering
$$1 \to F \to \tilde G \times T' \stackrel{\pi}{\to} G \to 1 \leqno{(12)}$$
defined over $k$.
Since the innner twist does not effect the
center it is obvious that we have also a special covering
$$ 1 \to F \to \tilde G_1  \times T' \to G_1 \to 1, \leqno{(13)}$$
where $\tilde G_1$ is the simply connected covering of the semisimple part
of $G_1$. The exact sequences (12) and (13) induce the following exact
sequences of groups of R-equivalences
$$  (\tilde G(k_v) \times T'(k_v))/R \to G(k_v)/R \to \H^1(k_v,F)/R \to 0,
\leqno{(14)}$$
$$ (\tilde G_1(k_v) \times T'(k_v))/R \to G_1(k_v)/R \to \H^1(k_v,F)/R \to 0,
\leqno{(15)}$$
(compare with (8) and (9)). 
Since the first groups in the exact sequences (14), (15)
are trivial, we obtain
$$G(k_v)/R \simeq G_1(k_v)/R \simeq \H^1(k_v,F)/R. $$
\\
Now assume that $G_1$ is an inner form of a group $G$ satisfying
3) $a)$ above. We will show that $G_1$ has weak approximation with respect to any
finite set $S$ of valuations $v$ where $G$ has metacyclic splitting field
extension $l_v/k_v$. By Proposition 3.4, it suffices to show that $G_1(k_v)/Br$
is trivial for all $v \in S$.
\\
Let $G = G'S$, where $G'$ is the semisimple part of $G$, and $S$ a central
torus. Let $F = S \cap G'$, $\tilde G$ be the simply connected
covering of $G'$, $\pi_G : \tilde G \to G'$
the canonical isogeny. As in part $a)$ we denote
$$F_1 = \{(f,f) : f \in F\} \hookrightarrow S \times G',$$
so we have a central isogenies
$$ 1 \to F_1 \to S \times G' \stackrel{\alpha}{\to} G=SG' \to 1,$$
$$ 1 \to F_2 \to S \times \tilde G \stackrel{\beta}{\to} S \times G' \to 1,$$
where $F_2 = \{(x,1) : x \in \Ker (\pi_G)\} \simeq \Ker \pi_G$.
Denote by $$\pi : 
S \times \tilde G \to SG'$$
 the composite of isogenies $\alpha$ and $\beta$. Then 
one checks
that
\begin{tabbing}
\hspace*{3cm}$\mu := \Ker \pi$ \=$= \{(g,\pi_G(g)) : g \in \pi_G^{-1} (F)\}$
\\
\\
\>$\simeq \pi_G^{-1} (F)$\\
\\
\>$\hookrightarrow \tilde F :=Cent(\tilde G)$
\end{tabbing} 
Since $G_1$ is an inner twist of $G$, $G_1 = S G_1'$, where $G_1'$ is the 
semisimple
part of $G_1$, and $S \cap G_1' = S \cap G = F$. 
We define 
$$\mu_1 : = \{(x,x) : x \in \mu\},$$ 
$$ H = (Z \times (S \times \tilde G))/\mu_1,$$
$$H_1 = (Z \times(S \times \tilde G_1))/\mu_1.$$
Then from the construction it follows that
$$H = \tilde G P, H_1 = \tilde G_1 P,$$
where $P$ is the connected center of $H$ and $H_1$, and also
$$\tilde G \cap P = \tilde G_1 \cap P. \leqno{(16)}$$ 
Therefore we have 
$$ T_0 :=H/\tilde G \simeq T_1 :=H_1/\tilde G_1.  \leqno{(17)}$$
From the proof of the part $a)$ above (Proposition 4.7)
we see that $T_1(k_v)/Br$ is trivial
since $T_0(k_v)/Br$ ($= H(k_v)/Br = G(k_v)/Br$)
 is trivial (see the proof of 4.7 of part $a)$). Therefore  by
Corollary 3.6.3 and relations
(16), (17) $H_1(k_v)/Br$ (
which is isomorphic to $T_1(k_v)/Br)$
 is also trivial, thus so is $G_1(k_v)/Br$.
The proof of Theorem 4.6 is complete.
\halm 
\\
\\
Now we consider some other analogs related with R-equivalence relations.
We have the following extension of similar property of tori (see Theorem 1.1, 
(4)).
\\
\\
{\bf 4.9. Theorem.} {\it Let G be a connected linear algebraic group defined
over a \p -adic field $k_v$. Then the group of R-equivalence classes and
the group of Brauer equivalence classes coincide :}
$$G(k_v)/R = G(k_v)/Br.$$
\Pr 
As above, we may assume that $v$ is non-archimedean and $G$ is reductive.
\\
\\
{\it Step 1. Let G be a connected reductive $k_v$-group, $G_q$ be its quasi-split
inner form defined over $k_v$. Then we have}
$$G(k_v)/R \simeq G_q(k_v)/R.\leqno{(18)}$$
This has been proved in Theorem 4.6, 3).
\\
\\
{\it Step 2. Let $G_q$ be a connected reductive quasi-split $k_v$-group. Then}
$$G_q(k_v)/R = G_q(k_v)/Br.\leqno{(19)}$$
\Pr
Take a maximal $k_v$-torus $T$ of $G_q$ containing a maximal $k_v$-split
torus $S$ of $G_q$. Then we have
$$T = Z_{G_q}(S),$$
and Bruhat decomposition for $G_q$ shows that (see [CTS1], Section 7)
we have
$$T(k_v)/R = G_q(k_v)/Br,$$
$$T(k_v)/Br = G_q(k_v)/Br.$$
Since for tori $T$ we have
$$T(k_v)/R = T(k_v)/Br,$$
by Theorem 1.1, 4), hence $G_q(k_v)/R = G_q(k_v)/Br.$
\\
\\
{\it Step 3. If $G_q$ is a quasi-split inner form of a connected reductive
$k_v$-group G, then }
$$G(k_v)/Br \simeq G_q(k_v)/Br.\leqno{(20)}$$
\\
Indeed, by Proposition 3.6 we may asume that the semisimple part $\tilde G$
(resp. $\tilde G_q$) of $G$ (resp. $G_q$) is
simply connected. From the proof of Theorem 4.6, 3) b (see (16), (17)),
it follows that
we have the following canonical isomorphism 
$$G_q/\tilde G_q \simeq G/\tilde G \simeq T,$$
where $T$ is quotient torus of $G$ (and $T$ is defined on the same field as $G$).
We know by Corollary 3.6.3 that
$$G_q(k_v)/Br \simeq T(k_v)/Br \simeq G(k_v)/Br,$$
hence (20) holds.
\\
\\
Now the theorem follows from the combination of (18), (19), (20). \halm 
\\
\\
{\bf 4.10. Remarks.} $1)$ In [CTS2], Remarque
2.8.17, it was shown that for a given
smooth variety $X$ over a local \p -adic field and under some condition
$(H'_{\lambda})$
on the universal torsor under some torus, the Brauer and R-equivalence are
the same. Also, it is a very general method to obtain such kind
of results (e.g. one may obtain similar results for tori over \p -adic fields
(see [CTS2], Section 2, for details). 
\\
\\ 
$2)$ All results above tell us that if a connected reductive group $G$ over
a number field $k$ fails to have weak approximation over $k$, then for some valuation $v$
(which is necessarily non-archimedean), and the quasi-split inner form $G_q$ of $G$
(which is necessarily non-split), we have $G_q(k_v)/Br \ne 1$.
\\ 
\\
$3)$ Our assumption in Theorem 4.6 on the existence of metacyclic
extension of $k_v$ splitting $G$ has local character, so it is weaker than 
that of Sansuc [Sa], Corollaire 5.4, p. 34.
\\
The following {\it local-global} statement (or principle) would show that
our result is equivalent to that of Sansuc :
\\
\\
{\it A connected reductive group G defined over a number field k has a
metacyclic splitting field if and only if it is so over all  
completions $k_v$ of k.}\\
\\
$4)$ Equally it is natural (and important)
 to ask for which class $\cal{G}$ of finite groups
the following holds. We say that a finite Galois extension $k'/k$ is
a $\cal{G}$-extension if $Gal(k'/k) \in \cal{G}$. We require that
$\cal{G}$ be a kind of formation of groups,
i.e., it is closed with respect to the operations 
of taking subgroups, factor groups and finite direct product.
Then we ask $when$ the following
holds :
\\
\\
{\it A connected reductive group over a
number field k  has a $\cal{G}$-splitting field
if and only if it is so over all comletions $k_v$.}
\\
\\
Now we are able to formulate and prove a close analog of the exact 
sequence $(Br)$ for groups of R-equivalence classes.
\\
\\
{\bf 4.11. Theorem.} {\it Let $G$ be a connected linear algebraic
group defined over a number 
field k. Let H be a z-extension of the reductive part
of G, T be its torus quotient and S be the Neron - Severi
torus of T. Then in the following exact sequence}
$$1 \to \Ker \rho_G \to G(k)/R \stackrel{\rho_G}{\to} \prod_v G(k_v)/R \to 
\A(G) \to 1
\leqno{(R')}$$
{\it the subgroup $\Ker\rho_G$ has finite index $n_T= [T(k)/R : \III(S)]$,
(or the same, $Card(T(k)/Br)$)
in $G(k)/R$.}
{\it Moreover the following sequence is exact} 
$$ 1 \to \Ker \rho_G \to G(k)/R \to G(k)/Br \to 1,$$
{\it and the image of $G(k)/R$ in $\prod_v G(k_v)/R$, being isomorphic
to $G(k)/Br$, is also isomorphic to $T(k)/Br$.}
\\
\\
\Pr
The fact that the first sequence is exact follows from Proposition 3.2.
From Theorem 3.7 we have 
the following commutative diagram
with exact rows
\\
\\
\hspace*{2cm}$1 \to \Ker \rho_G \to \hspace{0.2cm}G(k)/R \stackrel{\rho_G}{\to}
\hspace{0.15cm} \prod_v G(k_v)/R \to \A(G) \to 1$
\\
\\
(21)
$\hspace*{2cm}\downarrow \lambda_G \hspace{1.5cm}\downarrow \lambda_G' \hspace{1.5cm} 
\downarrow \simeq \hspace{2cm} \downarrow =$ 
\\
\\
\hspace*{2.8cm}
$1 \hspace{0.7cm}\to \hspace{0.4cm} G(k)/Br \stackrel{\gamma_G}{\to}
 \prod_v G(k_v)/Br \to \A(G) \to 1.$ 
\\
\\
In the above diagram, the homomorphism $\lambda_G$ is induced from $\lambda_G'$
since we have the vertical
isomorphism "$\simeq$" due to Theorem 4.9, and it is clear
that 
$$\lambda_G' (\Ker (\rho_G)) \subset \Ker (\gamma_G).$$ 
Therefore it follows that
$$G(k)/R /\Ker(\rho_G) \simeq G(k)/Br $$
and the image of $G(k)/R$ in the product $\prod_v G(k_v)/R$ is isomorphic to
the group
$G(k)/Br \simeq T(k)/Br$ and has order equal to $n_T = [T(k)/R : \III(S)]$ by Proposition 2.6
and Theorem 3.7.
\halm 
\\
\\
From Theorem 4.11 it follows that to determine the structure of
$G(k)/R$ one needs to understand the structure of $\Ker \rho_G$, which
is given in the following theorem.
We also  derive the following analog of the exact sequence (R)
in Section 1 in the case
the number field $k$ is totally imaginary and also in many other cases,
namely if the semisimple part of $G$ contains no anisotropic
factors of (exceptional, trialitarian)
type $\D_4$ nor $\E_6$.
\\
\\
{\bf 4.12. Theorem.} 1) {\it Let G be a connected linear
algebraic group defined over number field k. Then we have the following 
commutative diagram, where all rows and columns are exact sequences }
\\
\\
\\
\hspace*{1.9cm}
$\matrix{&&\tilde G(k)/R&\stackrel{=}{\to}&\tilde G(k)/R&&&&\cr
&&&&&&&\cr
\vspace*{0.5cm}
         &&\Big \downarrow&&\Big\downarrow \pi_R&&&&\cr
\vspace*{0.5cm}
  1&\to&\Ker \rho_G&\to&G(k)/R&\stackrel{\lambda'_G}{\to}&G(k)/Br&\to&1 \cr
\vspace*{0.5cm}
  &&\Big \downarrow p &&\Big \downarrow q&&r \Big \downarrow \simeq &&\cr
\vspace*{0.5cm}
  1&\to&\III(S)&\to&T(k)/R&\stackrel{\lambda'_T}{\to}&T(k)/Br&\to&1 \cr
\vspace*{0.5cm}
&&\Big \downarrow &&\Big \downarrow &&&&\cr
\vspace*{0.5cm}
  &&1&&1&&&&\cr}$
\\
\\
{\it and $\tilde G \stackrel{\pi}{\to} G_s$ is the simply connected covering
of the semisimple part $G_s$ of G, T is the torus quotient of the reductive
part of a z-extension H of G, S the Neron - Severi torus of T 
and  r is an isomorphism.}
\\
\\
\\
2) {\it If the semisimple part of G  contains no anisotropic 
almost simple factors 
of types $\D_4$ (trialitarian) nor $\E_6$ then $\tilde G(k)/R =1$.
In particular, p and q are isomorphisms,
$G(k)/R \simeq T(k)/R$ 
and the following
exact equence (R') holds 
for G.}
$$ 1 \to \III(S) \to G(k)/R \stackrel{\rho_G}{\to}
 \prod_v G(k_v)/R \to \A(G) \to 1. 
\leqno{(R')}$$
{\it  In general, (R') holds for all connected linear algebraic groups
G if and only if all simply connected almost simple groups have trivial
group of R-equivalence classes. } 
\\
\\
3) {\it If k is totally imaginary number field, then p, q
are also isomorphisms and the 
exact sequence (R') holds for G.}
\\
\\
First we need the following results.
\\
\\
{\bf 4.13. Theorem.} ([CM], Thm. 4.3.) {\it Let G be an almost simple 
algebraic group 
of outer type $^2\A_n$ defined over a field k, $G(k) =\SU(\Phi,D),$ where 
$\Phi$ is
the associated hermitian form with respect to an involution J of second kind
over a division algebra D of center K. Let $\Sigma_J$ (resp. $\Sigma'_J$) 
be the group of elements which are J-symmetric
(resp. with J-symmetric reduced norm)
of D. Then}
$$G(k)/R \simeq \Sigma'_J/\Sigma_J.$$ 
{\bf 4.14. Theorem.} [H3] {\it 
Assume that $k$ is totally imaginary number field. Then
any simply connected semisimple $k$-group has
trivial (Galois) 1-cohomology,
 and anisotropic almost simple k-groups are of type
$\A_n$.}
\\
\\
{\bf 4.15. Proposition.} [T3] {\it Let H be a z-extension of
a connected reductive group G, all defined over a field k. Then the natural 
projection $H \to G$ induces
 canonical isomorphism of abstract groups}
$$H(k)/R \simeq G(k)/R.$$ 
{\it Proof  of Theorem 4.12.}
\\
\\
1)
The assertion regarding the exactness
of two rows and that $r$ is an isomorphism in the above diagram
 follows from the commutative diagram (21), and the bijectivity 
of $r$ follows from
Theorem 3.7.
\\
\\
{\bf 4.16. Lemma.} {\it With above notation, $\Ker q \subset \Ker \rho_G$.}
\\
\\
\Pr~
Indeed, if $x  \in \Ker q$ then $r(\lambda'_G(x)) = \lambda'_T(q(x)) = 0,$
hence $\lambda'_G(x) =0$ since $r$ is an isomorphism. As we mentioned above,
the rows in the diagram are exact, so $x \in \Ker \rho_G$. \halm
\\
\\
The following is a well-known (and trivial) result from homological
algebra.
\\
\\
{\bf 4.17. Lemma.} {\it In the above diagram, p is surjective if and only if
q is surjective.}
\\
\\
{\bf 4.18. Lemma.} {\it Let T be a torus defined over a field k, S a finite
set of discrete valuations of k. Then}
$$Cl_S(RT(k)) = RT_S.$$
{\it In particular, if S consists of real valuations then $RT(k)$ is
dense (in the S-adic topology) in $T(k)$.}
\\
\\
\Pr~
Let
$$ 1 \to N \to P \stackrel{q}{\to} T \to 1$$
be a flasque resolution of $T$ over $k$. Here $P$ is an induced $k$-torus
and $S$ is flasque.
By [CTS1], Th\'eor\`eme 2, in the above exact sequence we have
$$q(P(k)) = RT(k),$$
hence 
$$ Cl_S(RT(k)) = Cl_S(q(P(k))).$$
If $x \in q(Cl_S(P(k)))$, $x = q(y),$ where $y = lim_n~ p_n,$ $p_n \in P(k)$
then 
\begin{tabbing}
\hspace*{3cm}$x$ \=$=q(lim_n~p_n)
=lim_n~q(p_n)$
\\
\\
\>$ \in  Cl_S(q(P(k))) = Cl_S(RT(k)),$ 
\end{tabbing}
hence $q(Cl_S(P(k))) \subset Cl_S(RT(k)).$
\\
Since $P$ has weak approximation property, $Cl_S(P(k)) = P_S,$
and as mentioned above, $q(P_S) = RT_S$, hence
$$RT_S \subset Cl_S(RT(k)).$$
On the other hand, by Proposition 2.1, $RT_S$ is an open subgroup of
$T_S$ containing $RT(k)$, hence the first assertion follows.
The rest of the lemma follows from the previous one and also
from the fact that any torus
over the real numbers are rational. \halm
\\
\\
{\bf 4.19. Lemma.} {\it With notation as in the theorem, q (hence
also p) is surjective.}
\\
\\
\Pr~
We need only show that 
$$T(k) = q(G(k)) RT(k),$$
where we may assume that $G$ has simply connected semisimple part $\tilde G$, 
and $q : G \to T = G/ \tilde G$ is the projection.
\\
By Lemma 4.18 we know that
$$T_{\infty} = Cl_{\infty}(RT(k)).$$
For $x \in T(k)$ we have
$$ x = lim_n~ r_n, r_n \in RT(k)$$
(the limit is taken with repsect to the archimedean $\infty$-topology).
We have the following commutative diagram similar to (6) :
\\
\\
\noindent
\hspace*{5cm}$G(k) \to T(k) \stackrel{\delta}{\to} \H^1(k,\tilde G)$
$$
\Big \downarrow \alpha \hspace{0.6cm} \Big
\downarrow \beta \hspace{0.8cm} \Big \downarrow \gamma
 \leqno{(22)}$$
\hspace*{4.8cm}
$G_\infty \hspace{0.2cm}\to T_\infty \stackrel{\delta_\infty}{\to} 
\prod_{v \in \infty} \H^1(k_v,\tilde G)$
\\
\\
where all rows are exact and all arrows are continuous with respect to
the topologies induced from $G_\infty$ and $T_\infty$. 
We have
\begin{tabbing}
\hspace*{3.5cm}$\delta_\infty(\beta (x))$ \=$=lim_n~ \delta_\infty(\beta(r_n))$
\\
\\
\>$=\delta_\infty(\beta (r_n)), \forall n >N,$
\end{tabbing}
for some fixed $N$, since $\prod_{v\in \infty} \H^1(k_v, \tilde G)$
is finite.
Hence 
\begin{tabbing}
\hspace*{3.5cm} $\delta_\infty(\beta (x))$ \=$= \gamma(\delta(x))$
\\
\\
\>$= \gamma(\delta(r_n))$
\end{tabbing}
for $n >N$. Since $\gamma$ is an isomorphism (i.e. bijection) by Hasse 
principle, one concludes that
$$\delta(x) = \delta(r_n), \forall n >N.$$
One checks, by using the interpretation of the coboundary map
$\delta$ (see [Se]) that $$x = r_n q(g),$$
for some $g \in G(k)$. Thus
$$ T(k) = RT(k) q(G(k),$$ i.e., $q$ 
(and $p$, by Lemma 4.17) is surjective.  \halm
\\
\\
With this lemma, the proof of the exactness of the first column in Theorem 4.12
is complete. Next we consider the exactness of the second column. We have
 the following general result.
\\
\\
{\bf 4.20. Lemma.} {\it Let k be a field of characteristic $0$, G
a connected reductive group with simply connected semisimple part $\tilde G$,
$T = G / \tilde G$. Then we have the following exact sequence of groups}
$$ \tilde G(k)/R \to G(k)/R \to T(k)/R.$$
\Pr~
It is obvious that if the lemma is true for some power $G^n = G \times \cdots
\times G,$ then it is also true for $G$, so by virtue of Lemma 4.8 we may
assume that $G$ has a special covering $\tilde G \times T'$, where $T'$
is an induced $k$-torus, and we have the following exacts sequence
of algebraic groups, all defined over $k$ :
$$ 1 \to F \to \tilde G \times T' \to G  \to 1,$$
where $F$ is a finite central subgroup of $\tilde G \times T'$. 
From this we derive the following $3\times 3$-commutative diagram
\\
\\
\hspace*{2cm}
$\matrix{&&1&&1&&1&&\cr
&&&&&&&\cr
 \vspace*{0.5cm}
         &&\Big \downarrow&&\Big \downarrow&&\Big \downarrow&&\cr
 \vspace*{0.5cm}
  1&\to&F \cap \tilde G&\to&F&\stackrel{u}{\to}&F'&\to&1 \cr 
 \vspace*{0.5cm}
  &&\Big \downarrow &&\Big \downarrow && \Big \downarrow &&\cr
 \vspace*{0.5cm}
  1&\to&\tilde G&\to&\tilde G \times T'&\to&T'&\to&1 \cr
 \vspace*{0.5cm}
 && l \Big \downarrow &&\Big \downarrow\pi&&\Big \downarrow&&\cr 
 \vspace*{0.5cm}
 1&\to&\tilde G&\to&G&\to&T&\to &1\cr
\vspace*{0.5cm}
&&\Big \downarrow && \Big \downarrow && \Big \downarrow && \cr
&&1&&1&&1&&\cr }$
\\
\\
Since $\tilde G$ is simply connected, $l$ is an isomorphism,
hence $F \cap \tilde G = 1$, and $u$ is also an isomorphism.
From the diagram above we derive the following commutative diagram
\\
\\
\hspace*{1cm}
$\matrix{(\tilde G(k)\times T'(k))/R&\stackrel{\pi_R}{\to}&G(k)/R&\to&
\H^1(k,F)/R\cr
&&&&\cr
\vspace*{0.5cm}
 \Big \downarrow p'&&\Big\downarrow q'&&\simeq \Big \downarrow r'\cr
 \vspace*{0.5cm}
 T'(k)/R&\to&T(k)/R&\to&\H^1(k,F')/R \cr }$
\\
\\
where all rows are exact (see Proposition 4.4.1), and $r'$ is an isomorphism.
Since $T'$ is an induced $k$-torus, $T'(k)/R =1$, and 
$$(\tilde G(k) \times T'(k))/R \simeq \tilde G(k)/R.$$
If $x \in G(k)/R$ such that $q'(x)=1$, then by chasing on this diagram
we see that $x \in \Im \pi$, thus 
$$\Ker q' = \Im (\tilde G(k)/R \to G(k)/R)$$
as required.  
\\
\\
From this the exactness of the second column of the diagram
in the theorem is proved, hence
we have finished the proof of 1).
\\
\\
\\
2) The "general" part of 2) follows directly from 1). Next we show that if
the semisimple part of $G$ does not contain anisotropic almost
simple factors of exceptional types $\D_4$ and $\E_6$ then (R') holds
for $G$.   
\\
 To see this, we reduce the proof
to the following situation. Namely, one may assume that the group
$G$ above is connected reductive with anisotropic
semisimple part.
Indeed, it is clear that we may assume $G$ to be reductive.
\\
\\
{\it Step 1.} {\it If H is an almost simple simply connected group
over k then $H(k)$ has no proper noncentral normal subgroups
except possibly for the following types:}
\\
$\bullet$ {\it anisotropic $\A_n$, $\D_4$ (exceptional), $\E_6$;}
\\
$\bullet$ {\it isotropic exceptional types $^2\E^{35}_6$, $^2\E_6^{29}$.} 
\\
\\
This is the result of many authors where the readers are refered to
Chapter 9, Section 9.1 of [PR] for further information. (Quite recently
Y. Segev and G. Seitz announced that the Platonov - Margulis conjecture
is true for the case $^1\A_n$.)
\\
\\
{\it Step 2.} {\it If H is as in Step 1, but H can be of anisotropic
type $\A_n$,  then $H(k)/R =1$.}
\\
\\
This follows from Step 1, the well-known fact that $RH(k)$
is an infinite normal subgroup of $H(k)$, Theorem 4.13 
in combination with results of Wang (that the group
$\SK_1(A) =1$ for any central simple algebra $A$ over number field $k$,
and the result of Platonov - Yanchevski\v i
(that $\Sigma'_J/\Sigma_J = 1$ for number field $k$).
\\
\\
{\bf 4.21. Lemma.} {\it If H is simply connected either of isotropic type
$^2\E_6^{35}$ or $^2\E_6^{29}$ then $H(k)/R =1$.}
\\
\\
\Pr
 If $S'$ is a maximal $k$-split torus of $H$, then
it is well-known by [CTS1] that we have the following functorial
isomorphism  of abstract groups
$$ H(K)/R \simeq Z_H(S')(K)/R,$$
for any field extension $K$ of $k$. Hence by Theorem 3.4, we have 
$$\A(H) \simeq \A(Z_H(S')). $$ 
(One can show in general that this last
isomorphism holds for any field $k$, see [T4].)
\\
{\it Case $^2\E_6^{35}$.} Let $S'$ be a maximal $k$-split torus of $H$. The
Tits index of $H$ is as follows
\noindent
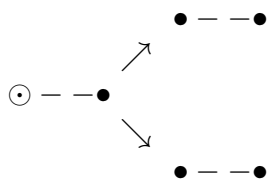
\begin{tabbing}\hspace*{6cm}\=$\bullet -- \bullet$
\\
\hspace*{5.3cm}$\nearrow$
\\
\hspace*{3.8cm}$\odot -- \bullet$
\\
\hspace{5.3cm}$\searrow$
\\
\>$\bullet -- \bullet$
\end{tabbing}
One can check that the centralizer $Z : =Z_H(S')$ of $S'$ in $H$ is 
$$ Z_H(S') = S' L,$$
where $L$ is an almost simple simply connected $k$-group of type $^2\A_5$.
Since $L(k)/R =1$ by Step 2, from the result of 1) (namely from the
commutative diagram in the theorem),
 we have  the following exact sequence (R') for
$Z$ :
$$ 1 \to \III(S) \to Z(k)/R \to \prod_v Z(k_v)/R \to \A(Z) \to 1, \leqno{(23)}$$
where $S$ is the Neron - Severi torus of the torus $Z/L$
since $Z$ is the $z$-extension of itself. Since $Z/L$ is a $k$-split torus, 
$S$ has trivial cohomology, so $\III(S)$ is trivial. As we notice
earlier that 
$$Z(k_v)/R \simeq H(k_v)/R$$
is trivial for all $v$ since $H$ is simply connected, hence
from (23) it follows that $$1 = Z(k)/R = H(k)/R$$
as required.
\\
{\it Case $^2\E_6^{29}$.} The Tits index of $H$ is as follows
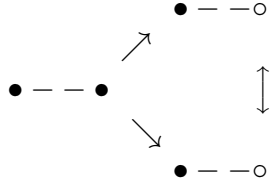
\begin{tabbing}
\hspace*{6cm}\=$\bullet -- \circ$
\\
\hspace*{5.3cm}$ \nearrow$
\\
\hspace*{3.8cm}$\bullet -- \bullet \hspace{1.83cm}\Big \updownarrow$
\\
\hspace*{5.3cm} $\searrow$
\\
\>$\bullet -- \circ$
\end{tabbing}
We have
$$Z : =Z_H(S') = S'T_0 L,$$
where $L$ is a simply connected $k$-group of (classical) type $\D_4$
(hence satisfies $L(k)/R =1$),
and $T_0$ is a one-dimesnional $k$-torus.
As above we have the exact sequence (R') for
the group $Z$. In this case, the torus quotient
$$ T = Z/L = S'T_0 / (S'T_0 \cap L)$$
is a two-dimensional $k$-torus, which is rational over $k$  by
a classical result of Voskresenski\v i [V1]. Therefore
$$T(k)/R = 1.$$  
By [CTS1], Proposition 19(ii), we have the following 
exact sequence
$$ 0 \to \III(T)~\tilde{} \to \Br_a X \to \prod_v \Br_a X_v \to T(k)/R~ \tilde
{} \to \III(S)~ \tilde {} \to 0.$$
(Here $X$ denotes a smooth compactification of $T$ over $k$ and
$\hat S = \Pic(\bar X)$.) In particular, $\III(S) =0$. Further we argue
as above to obtain that $Z(k)/R $ is trivial, hence so is $H(k)/R$. 
\halm
\\
\\
So from 1), Steps 1, 2 and from Lemma 4.21
it follows that the exact sequence (R') holds for $G$
except possibly the case the semisimple part of $G$ contains
anisotropic almost simple factors of exceptional types
$\D_4$ and/or $\E_6$. Hence 2) is proved.
\\
\\
3) We claim that
$H(k)/R =1$  for any almost simple simply connected group $H$ over
$k$. If $H$ is anisotropic then by Theorem 4.14, $H$ is of type $\A_n$
and the claim follows from Steps 1, 2 in 2). Also from there,
by combining with Lemma 4.21, we know that the
claim holds for any isotropic group $H$.
So over totally imaginary fields $k$, $H(k)/R$ 
is trivial for all simply connected semisimple $k$-groups $H$. 
To obtain (R') we may use the result of 1) and Theorem 4.11. We supply also
a proof of this fact independent of 1) as follows.
\\
By using the same argument as in the proof of Theorem 4.4 (or 4.6)
and by using the Harder's result on
the triviality of the Galois cohomology
of simly connected groups (Theorem 4.14), we can show that 
$$G(k)/R \simeq G_q(k)/R,$$
where $G_q$ is a quasi-split inner form of $G$ over $k$ and we may
assume also that $G$ has simply connected reductive part and that $G$ is
reductive. Let $T_q$ be  a 
maximal $k$-torus of $G_q$ containing a maximal $k$-split torus $S$ of
$G_q$.  By the same
argument as in the proof of 4.6, 4.7,  
it follows that $G$ and $G_q$ have isomorphic torus quotients $T$.
It is also well-known that for simply connected quasi-split semisimple
groups $G_q'$, any maximal torus containing a maximal $k$-split torus is
also quasi-split (i.e. induced) torus. Denote such a torus by $T_q'$.
Then we have
$$ 1 \to T_q' \to T_q \to T \to 1$$
is an exact sequence of $k$-tori, and since $T_q'$ is cohomologically
trivial. Since this is a $z$-extension, from Proposition 4.15 it follows that
$$T_q(k)/R \simeq T(k)/R,$$   
hence $$G(k)/R \simeq T(k)/R,$$
and this is true for any field extension of $k$. In particular, 
$$G(k_v)/R \simeq T(k_v)/R$$
for all $v$, hence from (R) we deduce the exact sequence
$$ 1 \to \III(S) \to G(k)/R \to \prod_v G(k_v)/R \to \A(G) \to 1, \leqno{(R')}$$
where $T$ is the torus quotient of any $z$-extension of $G$ and
$S$ is its Neron - Severi torus. 
\\
\\
The proof of Theorem 4.12 is therefore  complete.
\halm
\\ 
\\
We derive the following consequence from Theorem 4.12 which gives
a relation between $RG(k)$ and $RG_S$, extending
the corresponding result for tori (Lemma 4.18)  over number fields.\\
\\
{\bf 4.22. Theorem.} {\it Let k be a number field, S a finite set
of valuations of k, G a connected linear algebraic k-group. Then}
$$RG_S = Cl_S(RG(k)).$$
\Pr
We claim that
$$RG_S =  \pi(\tilde G_S) Cl_S(RG(k)),$$
{\it where $\pi  : \tilde G \to G_s$ is the universal covering of
the semisimple part $G_s$ of G.}
\\
\\
We may assume that $G$ is reductive.
We first show that $RG_S$ contains the set on right hanside.
It is known and easy to see  that $RG_S$ is an open subgroup in
$\prod_{v \in S} G(k_v)$, hence also closed. We know that for
simply connected groups $$\tilde G_S = R\tilde G_S = Cl_S(\tilde G(k))$$
hence $$\pi(\tilde G_S) = \pi(R\tilde G_S)  
\subset RG_S,$$
and 
$$\pi(\tilde G_S) Cl_S(RG(k)) \subset RG_S.$$ 
To prove the other inclusion, first we assume that $\tilde G$ is
the semisimple part of $G$.
Let $T = G/\tilde G$.
We have the following commutative diagram
\\
\\
\noindent
\hspace*{5cm}$G(k)\stackrel{p}{\to}
 T(k) \stackrel{\delta}{\to} \H^1(k,\tilde G)$
$$
\Big \downarrow \alpha \hspace{0.6cm} \Big \downarrow \beta \hspace{1cm} 
\Big \downarrow \gamma
 \eqno{}$$
\hspace*{4.8cm}
$G_S \hspace{0.2cm}\stackrel{p_S}{\to} T_S \stackrel{\delta}{\to} 
\prod_{v \in S} 
\H^1(k_v,\tilde G)$
\\
\\
Let $x \in RG_S$. Then $p(x) \in RT_S = Cl_S(RT(k))$ by 
Lemma 4.18, so $$p(x) = lim_n~r_n, ~r_n \in RT(k).$$
Then
$$\delta(r_n) \to \delta (p(x)) = 1,~ n \to \infty$$
hence $$\delta(r_n) = 1, \forall n >N,$$
for some fixed $N$. Therefore $r_n \in p(G(k))$,
$r_n = p(g_n), g_n \in G(k)$ for $n >N$.
Thus $$lim_n~ p(x g_n^{-1}) = 1.$$
Moreover, from Lemma 4.20 and the fact that
$p(g_n) \in RT(k)$ we deduce
$$g_n \in \tilde G(k) RG(k), ~ \forall n >N. \leqno{(24)}$$
Let $$G = \tilde G T',$$
where $T'$ is a $k$-torus. The natural isogeny
$$p' : \tilde G \times T' \stackrel{}{\to} G$$
induces an open map
$$p' : \tilde G_S \times T'_S \to G_S.$$
In particular, $\tilde G_S T'_S$ is an open subgroup
of $G_S$. Since $Cl_S(RT'(k)) =RT'_S$ by Lemma 4.18, and
$RT'_S$ is open in $T'_S$ by Proposition 2.1, it follows that
$Cl_S(RT'(k))$ is open in $T'_S$. Hence
$\tilde G_S Cl_S(RT'(k))$ is an open subgroup of $G_S$, and 
so is $\tilde G_S Cl_S(RG(k)).$ Let $V_n$, $n=1,2,...$ be a nested system
of open neighbourhoods of 1 in $T_S$ such that
$V_{n+1} \subset V_{n} $ for all $n$, 
$$\cap_n V_n =\{1\},$$ 
and
$p(xg_n^{-1}) \in V_n,$ for all $n.$
Then $$xg_n^{-1} \in p^{-1}(V_n), \forall n,$$
or for all $n$ we have
\begin{tabbing}
\hspace*{2.5cm}
$x \in p^{-1}(V_n)g_n$ \=$\subset p^{-1}(V_n) \tilde G(k) RG(k), $
\\
\\
\>$\subset p^{-1}(V_n) \tilde G_S Cl_S(RG(k))$
\end{tabbing}
by (24). Since
\begin{tabbing}
\hspace*{2.5cm}$\cap_n p^{-1}(V_n)$ \=$=p^{-1}(\cap_n V_n)$
\\
\\
\>$=p^{-1}(1)$
\\
\\
\>$= \tilde G_S,$
\end{tabbing}
so $(V_n)$ form a nested system of open neighbourhoods
of $\tilde G_S$. 
Therefore for some $N$,
$$p^{-1}(V_n) \subset \tilde G_S Cl_S(RG(k)), \forall n > N,$$
since $\tilde G_S Cl_S(RG(k))$ is an open subgroup of $G_S$.
Thus
$$x \in \tilde G_S Cl_S(RG(k)) = Cl_S(\tilde G(k)) Cl_S(RG(k)).$$
since, by Kneser's and Harder's
results,  $\tilde G$ has weak approximation property over $k$.
\\
In general  case, let $H$ be a $z$-extension of $G$,
$$ 1 \to Z \to H \stackrel{\pi}{\to} G \to 1,$$
where $\pi$ induces the covering isogeny $\tilde G \to G_s \subset G$.
By Proposition 4.15 the projection $\pi$
induces  surjections
$$ RH_S \to RG_S,~ RH(k) \to RG(k),$$
hence 
\begin{tabbing}
\hspace*{2cm}$RG_S$ \=$= \pi(RH_S)$
\\
\\
\>$= \pi(\tilde G_S) Cl_S(RH(k)))$
\\
\\
\>$=\pi(\tilde G_S) \pi(Cl_S(RH(k)))$
\\
\\
\>$\subset \pi(\tilde G_S) Cl_S(\pi(RH(k)))$
\\
\\
\>$= \pi(\tilde G_S) Cl_S(RG(k)),$
\end{tabbing} 
and the claim is proved. Since $R\tilde G(k)$ is
an infinite normal subgroup of $\tilde  G(k)$,
$Cl_S(R\tilde G(k))$ is an infinite normal subgroup of
$Cl_S(\tilde G(k)) = \tilde G_S$. It is known that $\tilde G_S$ has
no proper infinite normal subgroup (consequence of
Kneser - Tits conjecture over local fields). Thus
$$Cl_S(R\tilde G(k)) = \tilde G_S,$$
hence $$\tilde G_S \subset Cl_S(RH(k))$$
and $RG_S = Cl_S(RG(k))$. Theorem 4.22 is proved.  \halm
\\
\\
We derive also the following finiteness result of the group
of R-equivalence classes.
\\
\\
{\bf 4.23. Corollary.} [G2] {\it If k is a number field then
$G(k)/R$ is finite.}
\\
\\
\Pr
From well-known theorem of Margulis - Prasad it follows that if $\tilde G$
is a simply connected semisimple $k$-group then $\tilde G(k)/R$
is finite. From the commutative diagram in Theorem 4.12 
and from the finiteness of $T(k)/R$ (see [CTS1, Corol. 2, p. 200]
it follows that $G(k)/R$ is finite. \halm
\\
\\
{\bf 4.24. Corollary.} {\it Let G be an adjoint semisimple group
defined over a number field k and $\tilde G$ be its simply connected
covering. If $\tilde G(k)/R =1$ then $G(k)/R =1$. In particular, it is so
if G contains no anisotropic factors of exceptional
types $\D_4, \E_6$.}
\\
\\
\Pr
Let $F = \Ker (\tilde G \to G),$ $\tilde G_q, G_q$ be quasi-split inner
forms of $\tilde G, G$ respectively. Denote by $\tilde T_q$, $T_q$
their corresponding maximal $k$-torus containing maximally $k$-split
torus, which are known to be induced tori. From the exact sequence
$$ \tilde G(k)/R \to G(k)/R \to \H^1(k,F)/R,$$
and from  the assumption it follows that $G(k)/R \to \H^1(k,F)/R$
is injective. By considering the corresponding sequence for tori
$$ \tilde T_q(k)/R \to T_q(k)/R \to \H^1(k,F)/R \to 0,$$
(we use the fact that $\H^1(k,\tilde T_q) =0$) and from the fact that
these tori are induced, so we have trivially $\H^1(k,F)/R =0$, thus
$G(k)/R $ is trivial.  \halm
\\
\def\gK{{G(K)/R}}
\section {An application to norm principle}
Let $A$ be a central simple associative algebra of finite 
dimension over its center
$k$, we denote by $\K_1(A) = A^*/[A^*,A^*]$
the Whitehead group of $A$ and by 
$\SK_1(A) = \SL_1(A)/[A^*,A^*]$ the reduced Whitehead group of $A$.
Then by a theorem of Draxl [D], for any finite extension $K$ of $k$ we have 
a functorial
norm maps $$\K_1(A\otimes K) \to \K_1(A),~\SK_1(A \otimes K) \to \SK_1(A).$$
If $G$ denotes  the
algebraic $k$-group defined by $G(k) = \SL_1(A)$, then 
by a theorem of Voskresenski (see e.g. [V1]) $$ \SK_1(A) \simeq G(k)/R$$
so we have also a functorial norm map 
$G(K)/R \to G(k)/R.$
\\
If $A$ is endowed with an involution of the second kind $J$, 
which is not trivial on $k$,
then one can introduce the unitary Whitehead group 
$$\UK_1(A) = A^*/\Sigma_J,$$
and the reduced unitary Whitehead group
$$\USK_1(A) = \Sigma'_J/\Sigma_J,$$
where $\Sigma'_J$ is the subgroup of $A^*$ generated by elements with 
$J$-symmetric reduced norm and $\Sigma_J$ denotes the subgroup of $A^*$ generated by
$J$-symmetric elements. It is easy to show that in this case we also
have norm maps
$$\UK_1(A \otimes K) \to \UK_1(A),~ \USK_1(A \otimes K) \to \USK_1(A).$$
If $A$ is equipped with an isotropic nondegenerate $J$-hermitian form $\Phi$
and $k_0$ is the subfield of $J$-symmetric elements of $K$
then a result of  Yanchevski\v i 
(see e.g. [MY]) showed that for the group $G$ defined by
$G(k_0) = \SU(A,\Phi)$, we have
$$\USK_1(A) \simeq G(k_0)/R,$$
so we have also a norm map for $G(k_0)/R$ in this case.
\\
In general, as suggested by these examples, one may ask
\\
\\
{\it Is there a functorial norm map}
$$N_{k'/k} : G(k')/R \to G(k)/R,$$
{\it for any reductive group $G$ defined over a field $k$
and for any  finite extension $k'/k$ ?}
\\
\\
 For an arbitrary field $k$ it seems
to be a difficult question. Here we apply our results in previous
sections to answer this question in the case $k$ is a local or global field
of characteristic 0.
\\
It turns out that the investigation 
(and the existence) of the Norm Principle
for the group of R-equivalence classes touches upon deep and interesting 
properties of semisimple algebraic groups over fields and it has
also some interesting applications (see e.g. [ChM]).
\\
One of our main ingredients
is  the following
extention of a result of P. Deligne [De] about the existence of a norm
maps between certains factor groups of connected reductive groups over 
local and global fields. (In fact this result of Deligne motivates our
study in the norm principle. There is a more general result in [T1] but
we restrict ourselves only to the case which we need.) 
\\
\\
{\bf 5.1. Theorem.} [T1] {\it 
Let k be a
local or global field of characteristic 0.}
\\
\\
1) 
{\it Let $\pi: G \to H$ be a homomorphism of connected reductive
groups, all defined over  
 k, such that the restriction
of $\pi$ to the semisimple part $G'$ of G is an isogeny onto that of 
$H$. Then for any
finite extension $k'$ of k there exists a functorial corestriction (norm) map}
$$H(k')/\pi (G(k')) \to H(k)/\pi (G(k)).$$
2) {\it Let $1 \to F \to G \to H \to 1$  be an exact sequence of 
linear algebraic k-groups, with G and H  connected reductive groups and
F a central subgroup of G. Let $k'$ be any finite extension of k,
and $\delta_{k'}$, $\delta_k$ be boundary homomorphism}
$$H(k') \to \H^1(k',F),~ H(k) \to \H^1(k,F),$$
{\it respectively. Then the corestriction homomorphism}
$$Cor_{k'/k} :
\H^1(k',F) \to \H^1(k,F)$$
{\it  induces a functorial corestriction (norm) map}
$$N_{k'/k} : \delta_{k'}(H(k')) \to \delta_k (H(k)).$$
\noindent
{\bf 5.2. Remark.}
The proof of this result in
[T1] is based on cohomological consideration, and
the map we constructed comes indeed from  a {\it corestriction map} in certain
cohomology (in terminology of [De], from {\it Picard category}). It was also
the starting point for {\it abelianized Galois cohomology} of Borovoi
(to appear).
\\
\\	
The following result provides one application of Theorem 5.1.
\\
\\
{\bf 5.3. Theorem.} {\it Let G be a connected reductive group over a
local or global field k, $k'$ a finite extension of k. Assume 
that the $K$-rational points of the simply connected 
covering $\tilde G$ of the semisimple part $G_s$ of G
have images lying in $RG(K)$
when $K=k'$ or k.
Then there is a 
functorial (corestriction) norm map}
$$ G(k')/R \to G(k)/R.$$
{\it Proof.}
By Lemma 4.8
there exists induced  $k$-tori $P$, $Q$ and a natural number $n$
such that we have a central isogeny of connected reductive groups
$$ 1 \to F \to \tilde G^n \times P \stackrel{\pi}{\to} G^n \times Q \to 1,$$
where $F$ is a  finite central subgroup.
We denote by $G_1 = \tilde G \times P$, $G_2 = G^n \times Q$.
It is clear that it suffices to prove the theorem for $G_2$.
By assumption, we have 
$$\pi(G_1(K)) \subset RG_2(K)$$
for $K= k, k'$.
We have the following
induced exact sequence of cohomology
$$ G_1(K) \to G_2(K) \stackrel{\delta}{\to} \H^1(K,F),$$
and the following exact sequence (see Proposition 4.4.1)
$$ G_1(K)/R \to G_2(K)/R \to \Im~ \delta/R \to 1,$$
where $K$ is any field extension of $k$ and $\Im \delta/R$ denotes the
set of $R$-equivalence classes of $\Im(\delta)$
with the natural $R$-equivalence induced
from that of $G(K)$.
Now we take $K=k'$ and
$K =k$.
By taking the induced cohomology we have the following
commutative diagram with exact sequences rows,
where $R :=RG_2(k)$, $R':=RG_2(k')$
are subgroups in $G_2(k)$, $G_2(k')$, respectively,
which are generated by elements $R$-equivalent to 1.
\noindent
\\
\\
\hspace*{3cm}$1 \to R'/\Im~ \pi \to G_2(k')/\Im~ \pi \to
G_2(k')/R'  \to 1$
\\
\\
$\hspace*{4.3cm}\Big \downarrow \alpha \hspace{2cm}\Big \downarrow \beta
\hspace{1.8cm} (?) \Big \downarrow$
\\
\\
$\hspace*{3cm} 1\to R/\Im~ \pi \to G_2(k)/\Im~ \pi \to G_2(k)/R \to 1.$
\\
\\
Here we have used the assumption that the image of $G_1(K)$ under
$\pi$ lies in the group $RG_2(K)$ for $K=k$ and $k'$ and $\beta$
is the functiorial norm map which  exists by Theorem  2.3
and the functorial norm
map  $\alpha$ exists by a result of Gille [G1], Prop. 3.3.2.
Therefore one obtains also a functorial norm map
$(?) : G_2(k')/R \to G_2(k)/R$
as required. \halm
\\
\\
From results of previous sections we derive the following.
\\
\\
{\bf 5.4. Proposition.} {\it Let k be a local \p -adic field and G a 
connected linear k-group.
Then for any finite extension $k'$ of k there is a functorial norm map}
$$ G(k')/R \to G(k)/R.$$
{\it Proof.} Let $\tilde G$ be the simply connected covering of the semisimple
part of $G$. It follows from above that as any simply connected semisimple
$k$-group, $\tilde G$
 has trivial group of R-equivalence classes. By applying
Theorem 5.3 there is
a functorial norm map for $G$ as indicated.
\halm
\\
\\
Now we assume that $k$ is a number field.
\\
\\
{\bf 5.5. Proposition.} {\it Let G be a connected linear algebraic 
group over a number field k.
If any almost simple factor of G is neither of anisotropic
exceptional type $\D_4$, nor $\E_6$,
then for any finite extension $k'$ of k there
is a functorial  norm map}
$$ G(k')/R \to G(k)/R.$$
{\it Proof.}
It follows from Section 4 that the simply connected covering of the semisimple
part of $G$ has trivial group of R-equivalence classes. Hence
there exists a functorial norm map as indicated by Theorem 5.3.  \halm
\\
\\
We can prove the analog of Theorem 5.1, 1)
for the group of R-equivalence classes,
 i.e., the Corestriction (Norm) Principle in its
relative form.
\\
\\
{\bf 5.6. Proposition.} {\it Let $\rho : G \to H$
be a homomorphism of connected reductive
groups, all defined over a local or global field k of characteristic
0, such that the restriction
of $\rho$ to the semisimple part of G is a central isogeny onto that of 
H. Denote
$M_G(k') = G(k')/R$ for any field extension k' of k,
 $\rho_R $ $($ resp. $\rho_R')$ the induced homomorphism
$M_G(k) \to M_H(k)$ $($ resp.  $M_G(k') \to M_H(k'))$. 
Then for any finite extension $k'$ of k there is a
functorial  norm map}
$$ M_H(k')/\Im(\rho'_R) \to M_H(k)/\Im(\rho_R).$$
\\
{\it Proof.} There are few ways to prove 
the theorem and we just indicate one of them.
Let $\tilde G$ be the simply connected covering of the semisimple part 
$G'$ of $ G$,
hence also of $H'$. Let $\alpha : \tilde G \to G$
and $\beta : \tilde G \to H$ be the 
corresponding homomorphisms. 
\\
\\
{\it Step 1.} We show that there exists a norm map
$$ M_G(k')/\Im(\alpha'_R) \to M_G(k)/\Im(\alpha_R),$$
(and also similar statement for the pair $(H, \beta)$).
We know that there exist a natural number $n$, induced $k$-tori $P$ and $Q$
such that there is a central isogeny
$$ 1 \to F \to \tilde G^n \times P \stackrel{\pi}{\to} G^n \times Q \to 1,$$
all defined over $k$. Here $F$ is a finite central subgroup of 
$\tilde G^n \times P$ of 
multiplicative
type.
\\
 By Theorem 5.1, the Corestriction Principle for the  image of 
$$ \delta_k : (G^n \times Q)(k) \to \H^1(k,F)$$
holds, i.e., for any finite extension $k'$  we have 
$$ N_{k'/k}(\Im(\delta_{k'}) ) \subset \Im(\delta_k).$$
Since the similar also holds for the image of the set of elements
$R$-equivalent to 1 by [G1], Prop. 3.3.2, it follows that for the induced
norm map
$$ N: \H^1(k',F)/R \to \H^1(k,F)/R,$$
$$ A := \Im~ ((G^n \times Q)(k')/R \to \H^1(k',F)/R),$$
and
$$B : =\Im~ ((G^n \times Q)(k)/R \to \H^1(k,F)/R),$$
we have 
$$N(A) \subset B.$$
Since $P,Q$ are induced tori, they are all rational over $k$, so
$\pi_R$ is just $\alpha_R ^n$, and we can finish the proof by using the exact sequence
$$ (\tilde G^n \times P)(K)/R \stackrel{\pi _R}{\to} (G^n \times Q)(K)/R \to
\Im~ \delta_K/R \to 0$$
for any field extension $K$ of $k$.
\\
\\
{\it Step 2.} 
From Step 1 
we have the following commutative diagram with exact rows 
\\
\\
\hspace*{2cm}
$\matrix{M_G(k')/\Im \alpha'_R & \stackrel{\rho'_R}{\to} & M_H(k')/\Im \beta'_R
& \to & M_H(k')/ \Im \rho'_R \cr
&&&& \cr
\Big \downarrow N_1&& \Big \downarrow N_2 && \Big \downarrow (?) \cr
&&&& \cr
M_G(k)/ \Im \alpha_R & \stackrel {\rho_R}{\to} & M_H(k)/ \Im \beta_R
& \to & M_H(k)/ \Im \rho_R \cr}$
\\
\\
\\
where the functorial norm maps $N_1, N_2$ exist by Step 1,
hence the same is true for (?) 
and the proposition follows.  \halm.
\section { Remarks, problems and conjectures}
{\bf 6.1.} From Theorem 4.12 it follows that over an arbitrary
field $k$, the finiteness of groups
 $G(k)/R$ for connected  reductive groups $G$
depends only on the finiteness of tori and simly connected
groups over $k$. It seems natural to state the following
\\
\\
{\bf Problem 1.} {\it Study the finiteness of $G(k)/R$ for simply connected
groups G over finitely generated $($over the prime field$)$ k.}
\\
\\
{\bf Problem 2.} {\it Same problem as above, but only for purely transcendental
extensions of {\bf Q} , ${\bf Q_p}$, {\bf  R}, {\bf C}.}
\\
\\
{\bf 6.2.} It is natural to make the following 
\\
\\
{\bf Conjecture 1.} {\it  The exact sequence $(R')$ holds for any
connected linear algebraic group $G$ over any number field k.}
\\
\\
Notice from above that
this conjecture is equivalent to the following
conjecture 
\\
\\
{\bf Conjecture 2.} {\it $G(k)/R$ is trivial for anisotropic
 almost simple simply connected
group G
of exceptional type $\D_4$ or $\E_6$ 
 over a number field k.}
\\
\\
As we have seen from above, the last conjecture is a
consequence of another stronger conjecture due to Platonov and Margulis.
(See [PR, Chapter 9] for more information.)
\\
\\
{\bf 6.3.} In our earlier preprint [T3], we propose another way to express an 
exact sequence
connecting groups of R-equivalence classes, weak approximation obstruction
A$(G)$ and the Tate - Shafarevich group of some finite Galois module. It is 
clear that
the above exact sequence (R') is, in a sense, more true
(or natural) analog of the initial
exactsequence (R) for tori established by Colliot- Th\'el\`ene and Sansuc.
\\
\\
{\bf 6.4.}
All main results of this paper remain true in the case of global 
function field,
except perhaps few exceptions regarding purely inseparable isogenies. 
This case will
be considered, together with other related
problems, in another paper under preparation. 
\\
\\
{\it Acknowledgements.} This paper owes much of its existence due
to works by Colliot - Th\'el\`ene and Sansuc.
The preliminary ideas and
first results of this paper were obtained while the author was enjoying a
Fellowship at McMaster University
(thanks to Professors I. Hambleton,
M. Kolster, C. Riehm, and  V. Snaith),
Lady Davis  Fellowship at Technion (thanks to Professor J. Sonn).
The first  version with further results of the paper was written
while the author was enjoying the stay at ICTP, Trieste, Italy (thanks to
Professor M. S. Narasimhan and Professor A. Kuku).
I would like to thank them and these
institutions  for the support and hospitality and all that helped to realize 
this work.
I would like to thank 
Professor J. -L. Colliot - Th\'el\`ene, Professor M. Borovoi,
Professsor B. Kunyavski\v i
and Professor J. Oesterl\'e for useful discussions related
to this paper.


\begin{thebibliography}{CTS1}
\bibitem[CM]-V. Chernousov and A. Merkurjev, R-equivalence and special unitary
groups, Preprint, Bielefeld Univ. No. 96083, 1996.
\bibitem[CTS1]-J. -L. Colliot-Th\'el\`ene et J. -J. Sansuc, La R-\'equivalence
sur les tores, Ann. Sc. E. N. S. t. 10 (1976), 175 - 230.
\bibitem[CTS2]-J. -L. Colliot-Th\'el\`ene et J. -J. Sansuc, La descente sur
les variet\'etes rationnelles, II, Duke Math. J. v. 54 (1987), 375 - 492.
\bibitem[De]-P. Deligne, Vari\'et\'es de Shimura : Interpretation
modulaire et techniques de construction de mod\`eles canoniques, Proc.
Sym. Pure Math. A.M.S. 33 (1979), Part 2, 247 - 289.
\bibitem[D]-P. Draxl, {\it Skew Fields},
 London. Math. Soc. Lecture Notes Ser. 1980.
\bibitem[G1]-P. Gille, R-\'equivalence et principe de norm
en cohomologie galoisienne, C. R. Acad. Sci. Paris t. 316, s\'er. I 
(1993), 315 - 320.
\bibitem[G2]-P. Gille, Une th\'eor\`eme de finitude arithm\'etique 
sur les groupes r\'eductifs,
C. R. Acad. Sci. Paris t. 316, s\'er. I (1993),  701 - 704.
\bibitem[H1]-G. Harder, Eine Bemerkung zum schwachen
Approximationsatz, Archive der Math. 19 (1968) 465 - 471.
\bibitem[H2]-G. Harder, Halbeinfache Gruppenschemata \"uber Dedekindringen,
Inv. Math. 4 (1967), 165 - 191.
\bibitem[H3]-G. Harder, \" Uber  die Galoiskohomologie der 
halbeinfacher Matrizengruppen,
I. Math. Z. 90 (1965), 404 - 428; II, ibid., 92 (1966), 396 - 415.
\bibitem[Kn1]-M. Kneser, Schwache Approximation in algebraischen Gruppen,
in {\it Colloque sur la th\'eorie des groupes alg\'ebriques},
 CBRM, Bruxelles, 1962, 41 - 52.
\bibitem[Kn2]-M. Kneser, Galois-Kohomologie halbeinfacher algebraischer Gruppen
\"uber \p -adischen K\"orpern, II, Math. Z. 89 (1965), 250 - 272.
\bibitem[M1]-Yu. I. Manin, Le groupe de Brauer - Grothendieck en g\'eometrie
diophantienne, Actes du Congr\`es intern. math. Nice 1 (1970) 401 - 411. 
\bibitem[M2]-Yu. I. Manin, {\it Cubic forms : Algebra,
Geometry, Arithmetic},
North Holland, 1986, 2nd ed. 
\bibitem[MT]-Yu. I. Manin and M. A. Tsfasman, Rational varieties : algebra, geometry and
arithmetic, Russian Math. Surveys v. 41 (1986), 51 - 116.
\bibitem[MY]-A. Monastyrnyi and V. Yanchevski\v i, Whitehead groups of 
spinor groups, Soviet Math. Dokl.v. 42 (1991), 351 - 355.
\bibitem[PR]-V. Platonov and A. Rapinchuk, {\it Algebraic 
Groups and Number Theory,} Academic Press, 1994.
\bibitem[O]-T. Ono, On relative Tamagawa numbers, Ann. Math. 82 (1965),
88 - 111.
\bibitem[S]-J. -J. Sansuc, Groupes de Brauer et arithm\'etique
des groupes alg\'ebriques
lin\'eaires sur un corps de nombres, J. reine. angew. Math. Bd. 
327 (1981), 12 - 80.
\bibitem[Se]-J. -P. Serre, {\it Cohomologie Galoisienne}, Lecture Notes in 
Math., v. 5, (Fifth edition), Springer - Verlag, 1994.
\bibitem[T1]-Nguyen Q. Thang, Corestriction Principle in Higher Dimensions,
I, II, Preprint, McMaster University, Sept. 1995, No. 2.
\bibitem[T2]-Nguyen Q. Thang, Rationality of almost simple algebraic groups,
Preprint, McMaster University, Sept. 1995, No. 8.
\bibitem[T3]-Nguyen Q. Thang, Weak approximation, R-equivalence and Whitehead 
group, in : {\it Algebraic K-Theory}, Fields Inst. Communications (A. M. S.)
v. 16 (1997), 345 - 354. (Cf. also a version of this paper in E-Preprint series,
Duke University, alg-geo:9603005, 1996.) 
\bibitem[T4]-Nguyen Q. Thang, On weak approximation 
in algebraic groups and varieties defined by a system of forms, J. Pure and App.
Algebra, v.113 (1996), 67 - 90.
\bibitem[Ti]-J. Tits, Classification of algebraic semisimple groups,
in : {\it Algebraic groups and Discontinuous groups}, 
Proc. Sym. Pure Math. A.M.S.
v. 9 (1966), 33 - 62.
\bibitem[V1]-V. E. Voskresenski\v i, {\it Algebraic tori}, Moscow, Nauka,
1977 (in russian).
\bibitem[V2]-V. E. Voskresenski\v i, Problems of R-equivalence in semisimple
groups, J. Soviet Math. 17 (1981), 1275 - 1287.
\end{thebibliography}
\end{document}